%% file: main.tex
\documentclass[11pt,a4paper,USenglish]{article}
\usepackage{template}
\usepackage{fix-cm}
\usepackage{document_style}
\usepackage{indentfirst}

\usepackage[maxbibnames=100]{biblatex}
\title{The Interplay of Shifted Square and Maximal Function Estimates \\ in the Context of Multilinear Fourier Multipliers}
\author{Andrew Haar\inst{1}}
\address{andrew.haar674@student.cuni.cz}

\begin{document}

\maketitle

\begin{abstract}
Following their appearance in 2014, so-called shifted square and maximal functions have seen an eruption of use in the study of singular integral operators. In this paper, we will generalize a recent theorem of G. Dosidis, B. Park, and L. Slav\'ikov\'a, which gave a sharp boundedness criterion for certain bilinear Fourier multipliers, to the general multilinear setting. In so doing, we will witness how the combined use of shifted square and maximal functions causes a loss of sharpness; we, then, repair this through a trick, which allows us to remove the shift from the square functions, placing it purely on the maximal functions. As an application to our main theorem, we establish the boundedness of certain singular integrals with rough homogeneous kernels lying in the Orlicz space $L(\log L)^\alpha$ when restricted to the unit sphere. This represents an edge case to what was previously known in the literature.
\end{abstract}

\import{Sections}{intro.tex}

\import{Sections}{setup.tex}

\import{Sections}{proof.tex}

\import{Sections}{counterexample.tex}



\printbibliography

\end{document}

%% file: Sections/intro.tex
\section{Introduction}\label{sec:intro}

One of the venerable problems facing the harmonic analysis community is the question of establishing necessary and sufficient conditions for the $L_p$-boundedness of Fourier multipliers. Unfortunately, it has been apparent for quite some time that we do not possess the tools to approach this problem. As such, the focus has shifted towards establishing (sharp) sufficient conditions, which is the object of this paper.

There are a number of branches in the lineage of multiplier theorems. As far back as 1960, L. H\"ormander \cite{hormander1960estimates} famously gave a criterion sufficient to guarantee the $L_p$-boundedness of a Fourier multiplier, which was based on the smoothness of the multiplier (along with appropriate decay assumptions for the derivatives). Later developments, initiated by Calder\'on and Torchinsky \cite{calderon1977parabolic}, loosened the requirement of classical smoothness to that of Sobolev and, further, brought the Hardy spaces into the picture. Through further efforts of many mathematicians, this progression of results has recently culminated in sharp (or almost sharp) sufficient conditions on the Sobolev smoothness of the multiplier, not just in the linear, but also in the multilinear setting, see \cite{grafakos2017hormander,grafakos2013multilinear,lee2021hormander,tomita2010hormander,} and the references therein.

While it is natural to attempt to place conditions on the multiplier in this pursuit, a related and fruitful strategy is to, instead, restrict the kernel, i.e. the inverse Fourier transform of the multiplier. Drawing from H\"ormander's criterion, the demand for smoothness of the multiplier in phase space translates to integrability criteria for our kernel. This was done by A. Baernstein and E. Sawyer \cite{baernstein1985embedding} for $H_1$, and was later extended by A. Seeger \cite{seeger1990remarks} to include the result for $L_p$. In fact, after turning their attention towards the kernel, rather than the multiplier, A. Baernstein and E. Sawyer observed \cite[Section 3]{baernstein1985embedding} that their results imply the classical ones of L. H\"ormander as well as that of Calder\'on and Torchinsky, which clearly establishes the merit of this approach.

The time has come for the boundedness criteria based on the kernel to catch up to those on the multiplier. That is, just as has been done for the H\"ormander condition, we wish to achieve sharp sufficient conditions in the full multilinear setting. As we shall see, actualizing this goal is made possible by the recent introduction of shifted variants of classical square and maximal functions.\footnote{The shifted square functions were first introduced by C. Muscalu in \cite{muscalu2014calderon1}, but the shifted maximal functions had appeared in the literature before at least as early as the 1980s, see \cite[Chapter 2, Section 5.10]{stein1993harmonic} and the references therein. After a small flurry of interest in them, however, they were forgotten until being revitalized in \cite{muscalu2014calderon1}.} Indeed, in a series of three papers \cite{muscalu2014calderon1,muscalu2014calderon2,muscalu2014calderon3}, C. Muscalu defined these shifted operators before proceeding in applying them to give a new proof for the $L_p$-boundedness of the Calder\'on commutators. Incidentally, it was an observation of R. Coifman and Y. Meyer \cite{coifman1975commutators} that Calder\'on's first commutator is more naturally viewed as a bilinear operator, which first sparked interest in multilinear harmonic analysis.

Since their debut in 2014, Muscalu's shifted square and maximal functions have seen a great deal of success past the specific case of the Calder\'on commutators, i.e. in establishing boundedness criteria for general multilinear operators. This brings us to the aforementioned imminent expansion of Baernstein and Sawyer's results to the multilinear setting: by taking inspiration from \cite{muscalu2014calderon1}, G. Dosidis, B. Park, and L. Slav\'ikov\'a \cite[Theorem 2]{dosidis2024boundedness} gave a new sharp criterion for the boundedness of certain fairly general bilinear Fourier multipliers by way of shifted square functions. The final piece of the puzzle is, thus, to do the same for multilinear operators, trilinear and beyond, which we shall do in this paper. While the proof for the bilinear case can be carried over to the $n$-linear setting, it fails to produce the sharp result that we want. For this reason, we have to apply a new technique, which allows us to avoid the use of shifted square functions in favor of the shifted maximal functions. We discuss this more precisely at the end of the introduction.

\begin{rem}
These shifted operators have been used in a wide array of results in harmonic analysis, especially in regard to singular integrals, see, e.g., \cite{dosidis2024multilinear,gaitan2020boundedness,lee2021hormander,sjogren1986admissible}. In particular, looking at \cite[Section 4]{lee2021hormander}, using a process inspired by \cite{muscalu2004bi}, which is similar in character to our methods below in Section \ref{sec:proofSetup}, the authors used estimates with shifted operators to establish their multilinear variant of H\"ormander's multiplier theorem.
\end{rem}

Now that we have oriented ourselves in the literature, it is time to state our main theorem. First, let us define some spaces; the spatial dimension, $d\in\nN$, will be fixed throughout the paper. As usual, $\cS(\nR^d)$ will be the space of rapidly decaying functions. Further, $\cS_0(\nR^d)\subset \cS(\nR^d)$ will be those Schwartz functions, $f$, all of whose moments vanish, i.e.
\begin{align*}
    \int_{\nR^d} x^\alpha f(x) \, dx = 0
\end{align*}
for all multiindices, $\alpha\in\nN_0^d$. It is well-known that $\cS_0(\nR^d)$ is dense in $H_p(\nR^d)$, the Hardy space, for $0<p<\infty$, see \cite[Section 5.1.5]{triebel1983function}. This will be our space of test functions, because $\cS$ is not dense in the Hardy spaces; it will only be a technical/cosmetic difference throughout the exposition that is to follow. We recall, as well, that $H_p(\nR^d) = L_p(\nR^d)$ if $1<p\leqslant \infty$. Finally, for brevity, we write
\begin{align*}
    Y_p(\nR^d) = \begin{cases}
        L_p(\nR^d) & 1\leqslant p<\infty \\
        BMO(\nR^d) & p=\infty .
    \end{cases}
\end{align*}
As usual, if $1\leqslant p\leqslant \infty$, then $p'$ is the H\"older conjugate of $p$, i.e. $\frac{1}{p} + \frac{1}{p'} = 1$.

Our object of interest, for $n\geqslant 2$ fixed, will be the following $n$-linear Fourier multiplier operator
\begin{align}
    T(f_1,\ldots,f_n)(x) := \int_{\nR^{nd}} m(\xi_1,\ldots,\xi_n) \left( \prod_{k=1}^n \widehat{f}_k(\xi_k) \right) e^{2\pi i (x,\xi_1 + \cdots + \xi_n)} \, d\xi , \quad x\in\nR^d , \label{eq:triop}
\end{align}
where $f_1,\ldots,f_n\in \cS_0(\nR^d)$. We will enforce that the symbol, $m$, is of the form
\begin{align}
    m = \sum_{j\in\nZ} m_0(2^j\cdot) , \label{eq:kernelStructure}
\end{align}
where $m_0$ is a bounded function with support in $\{(\xi_1,\ldots,\xi_n)\in \nR^{nd} : 1/2 \leqslant |(\xi_1,\ldots,\xi_n)| \leqslant 2\}$. With this setup, we observe that by writing $K = m_0^\vee$ and $K_\ell = 2^{n\ell d}K(2^\ell\cdot)$, the usual $L_1$ dilation, we have
\begin{align}
    T(f_1,\ldots,f_n)(x) = \sum_{\ell\in\nZ} K_\ell * (f_1 \otimes \cdots \otimes f_n)(x,\ldots,x) . \label{eq:multopcon}
\end{align}
The limitation in generality forced by our choice in structure for the multiplier is minor; we will say more on this while discussing Theorem \ref{thm:app}.

To state our theorem, we will also need the notation $\mx{k}(A)$ to mean the $k$-th largest number in a finite set, $A\subset\nR$;\footnote{We have a slight abuse of notation here, because the sets that we use this on can have repetition. For our purposes, duplicates will count independently. For example, if $A = \{3,3,2,1\}$, then $\mx{3}(A) = 2$ not $1$. The order, in which we count duplicates, will not matter for the argument, but we will be dealing with Cartesian points of the form $(\frac{1}{p_1},\ldots,\frac{1}{p_n},\frac{1}{p'})$, so we could, e.g., take the larger in a duplicate pair to be the one with the smallest coordinate, i.e. reading from left to right.} we define $\mn{k}(A)$ similarly.

\begin{thm}\label{thm:main}
Suppose $1 \leqslant p_1,p_2,\ldots,p_n,p \leqslant \infty$ satisfy
\begin{align*}
    \frac{1}{p} = \frac{1}{p_1} + \cdots + \frac{1}{p_n} .
\end{align*}
Further, take $K$ to be an integrable function in $\nR^{nd}$ with
\begin{align*}
    \supp \widehat{K} \subset \{(\xi_1,\ldots,\xi_n) \in \nR^{nd} : 1/2 \leqslant |(\xi_1,\ldots,\xi_n)| \leqslant 2 \} .
\end{align*}
Define
\begin{align*}
    D_\lambda(K) := \int_{\nR^{nd}} |K(y)| (\log(e + |y|))^\lambda \, dy
\end{align*}
and assume that this quantity is finite for $K$. The $n$-linear operator, $T$, is given through $K$ as in \eqref{eq:multopcon}. If
\begin{align}
    \lambda \geqslant \sum_{k=1}^{n-1} \mx{k} \left\{ \frac{1}{p_1} , \ldots , \frac{1}{p_n} , \frac{1}{p'} \right\} , \label{eq:lambdathm}
\end{align}
then there is a constant $C = C(p_1,\ldots,p_n,d,n) > 0$ such that
\begin{align}
    \|T(f_1,\ldots,f_n)\|_{Y_p(\nR^d)} \leqslant C D_\lambda(K) \prod_{k=1}^n \|f_k\|_{H_{p_k}(\nR^d)} \label{eq:maineq}
\end{align}
for all $f_1,\ldots,f_n \in \cS_0(\nR^d)$.

Conversely, if \eqref{eq:lambdathm} does not hold, then for every $C > 0$, there exist $f_1,\ldots,f_n\in\cS_0(\nR^d)$ along with an integrable kernel $K$ as above, such that
\begin{align*}
    \|T(f_1,\ldots,f_n)\|_{Y_p(\nR^d)} > C D_\lambda(K) \prod_{k=1}^n \|f_k\|_{H_{p_k}(\nR^d)} .
\end{align*}
\end{thm}

\begin{rem}
We remark that the theorem above is valid, as well, for $f_1,\ldots,f_n$ taken from $\cS(\nR^d)$ instead of $\cS_0(\nR^d)$ if we are willing to throw out the edge cases involving $H_1(\nR^d)$, i.e. if we assume $1<p_1,\ldots,p_n\leqslant\infty$ and $1\leqslant p<\infty$.
\end{rem}

One may wonder what happens in the linear case, i.e. when $n=1$. For this, there is an analogous result, which can be found in the same paper, which originally gave the bilinear result, check \cite[Proposition 1]{dosidis2024boundedness}, although, as we indicated earlier, it was first studied much earlier by A. Baernstein and E. Sawyer \cite[Theorem 3b]{baernstein1985embedding} when $p=1$ and later by A. Seeger \cite[Theorem 2.1]{seeger1990remarks} for all $1\leqslant p\leqslant \infty$. In both of these earlier papers, the authors did not additionally require the structure \eqref{eq:kernelStructure} of the kernel, which resulted in a slightly different condition on $\lambda$.

Briefly, the result in \cite[Proposition 1]{dosidis2024boundedness} is as one may expect from looking at the multilinear case above except for two details. First, when $p=1$, the more precise statement $T : H_1\to H_1$ is proven, rather than just $T : H_1 \to L_1$. Second, and most importantly, the sharp condition on $\lambda$ when $n=1$ is
\begin{align*}
    \lambda \geqslant \left|\frac12 - \frac{1}{p}\right| ,
\end{align*}
which comes from the shifted square function estimates, see Lemma \ref{lem:shiftedsquare}. When the structure on the kernel \eqref{eq:kernelStructure} was not required in \cite{baernstein1985embedding,seeger1990remarks}, the authors had to assume that $\lambda > \left|\frac12 - \frac{1}{p}\right|$, which is, in fact, sharp for the type of kernel they were using, see \cite[Theorem 4.2]{park2019fourier}.

This is in particular contrast to the multilinear case, because, at least when $p=2$, $\lambda = 0$ is permitted in the linear setting; in other words, the mere assumption that $K\in L_1$ is enough when $p=2$ and $n=1$. When $n > 1$, $\lambda$, restricted by \eqref{eq:lambdathm}, cannot be smaller than $\frac{n-1}{n+1}$, which means it never suffices to only assume $K\in L_1$.

In terms of applications, following the lead of \cite[Theorem 3]{dosidis2024boundedness}, our main theorem can be applied using a method of J. Duoandikoetxea and J. Rubio de Francia \cite{duoandikoetxea1986maximal} to obtain improved estimates for a well-known class of multilinear singular integral operators having rough homogeneous kernels. Broadly speaking, the process involves decomposing the rough kernel into a sum of smooth kernels, which accommodate the assumptions of Theorem \ref{thm:main}.

The rough kernel in question is built with $\Omega\in L_1(\nS^{nd-1})$, which has mean zero, i.e.
\begin{align*}
    \int_{\nS^{nd-1}} \Omega(\theta) \, d\nu(\theta) = 0 ,
\end{align*}
where the measure, $d\nu$, is the normalized surface measure on $\nS^{nd-1}$. Our $n$-linear operator is, then,
\begin{align}
    T_\Omega (f_1,\ldots,f_n)(x) = \text{pv}\, \int_{\nR^{nd}} \frac{\Omega(\theta)}{|y|^{nd}} f_1(x-y_1) \cdots f_n(x-y_n) \, dy , \quad x\in\nR^d , \label{eq:roughOperator}
\end{align}
where $\theta = y/|y|\in \nS^{nd-1}$ (and $y = (y_1,\ldots,y_n)\in \nR^{nd}$).

Boundedness properties for multilinear operators of the above form have been thoroughly studied in the last decade, dating back to work of L. Grafakos, D. He, and P. Honz\'ik \cite{grafakos2018rough} in the bilinear case. After a number of papers on bilinear operators \cite{grafakos20182,he2023improved}, researchers moved on to the multilinear setting as we have it here, see \cite{dosidis2024multilinear,grafakos2024multilinear}. The focus, so far, has been on boundedness properties of $T_\Omega$ with $\Omega\in L_q(\nS^{nd-1})$ for some $q > 1$; the accomplishment of the theorem below is to push this further to allow the edge case, where $\Omega$ is in certain Orlicz spaces, denoted $L(\log L)^\alpha(\nS^{nd-1})$; this is new in all but the bilinear case.

We say that $\Omega \in L(\log L)^\alpha(\nS^{nd-1})$ if the norm
\begin{align*}
    \|\Omega\|_{L(\log L)^\alpha(\nS^{nd-1})} = \inf \left\{ \lambda > 0 : \int_{\nS^{nd-1}} \frac{|\Omega(\theta)|}{\lambda} \left( \log\left( e + \frac{|\Omega(\theta)|}{\lambda} \right) \right)^\alpha \, d\nu(\theta)\leqslant 1 \right\}
\end{align*}
is finite. In particular, $\Omega \in L(\log L)^\alpha(\nS^{nd-1})$ if it satisfies
\begin{align}
    \int_{\nS^{nd-1}} |\Omega(\theta)| (\log(e + |\Omega(\theta)|))^\alpha \, d\nu(\theta) < \infty , \label{eq:omegaCondition}
\end{align}
which solicits the use of Theorem \ref{thm:main}. 

It is also worth noting that any $\Omega\in L_q(\nS^{nd-1})$ with $q > 1$ fulfills \eqref{eq:omegaCondition} as well, solidifying the Orlicz space as the alleged edge case. Pushing this further, as was similarly pointed out in \cite{dosidis2024boundedness}, if $\Omega\in L_q(\nS^{nd-1})$ for some $q > 1$, necessary and sufficient conditions have been established on $1<p_1,\ldots,p_n<\infty$ such that $T_\Omega : L_{p_1}\times\cdots\times L_{p_n}\to L_p$ is bounded, see \cite[Theorem 2 and Proposition 3]{dosidis2024multilinear}.\footnote{We are always assuming $\Omega$ has mean zero and that $\frac{1}{p} = \frac{1}{p_1} + \cdots + \frac{1}{p_n}$, although $p<1$ is permissible in the reference given.} We will not write out the full requirement, but it must hold that
\begin{align*}
    \frac{1}{p} + \frac{n-1}{q} < n .
\end{align*}
Therefore, if we are to prove a result for $\Omega\in L(\log L)^\alpha(\nS^{nd-1})$, then it will hold for all $\Omega\in L_q(\nS^{nd-1})$ with any $q > 1$, which implies that we must have $p\geqslant 1$. The result below contains the converse implication, at least when $\alpha \geqslant 2$. It, further, gives the result for a restricted range of points, $(\frac{1}{p_1},\ldots,\frac{1}{p_n})$, so long as $\alpha \in [\frac{2n}{n+1},2)$.

\begin{thm}\label{thm:app}
Let $1\leqslant p < \infty$ and $1<p_1,\ldots,p_n < \infty$ be such that
\begin{align*}
    \frac{1}{p} = \frac{1}{p_1} + \cdots + \frac{1}{p_n} .
\end{align*}
Suppose $\Omega\in L(\log L)^A(\nS^{nd-1})$ for some real number, $A$, with
\begin{align*}
    A \geqslant 1 + \sum_{k=1}^{n-1} \mx{k} \left\{ \frac{1}{p_1} , \ldots , \frac{1}{p_n} , \frac{1}{p'} \right\} = 2 - \sum_{k=1}^2 \mn{k} \left\{ \frac{1}{p_1} , \ldots , \frac{1}{p_n} , \frac{1}{p'} \right\} .
\end{align*}
We also require $\Omega$ to have mean zero. Then, for some $C = C(p_1,\ldots,p_n,d,n) > 0$, we have, for $T_\Omega$ as in \eqref{eq:roughOperator},
\begin{align*}
    \|T_\Omega (f_1,\ldots,f_n)\|_{L_p(\nR^d)} \leqslant C \|\Omega\|_{L(\log L)^A(\nS^{nd-1})} \prod_{k=1}^n \|f_k\|_{L_{p_k}(\nR^d)}
\end{align*}
for all $f_1,\ldots,f_n\in\cS(\nR^d)$.
\end{thm}

In this paper, we will refrain from giving a proof for the above theorem; instead, we encourage the reader to check \cite[Section 6]{dosidis2024boundedness}, where the bilinear case is established. The proof there has, other than minor aesthetic differences,\footnote{These differences mostly revolve around integrating over $\nR^{nd}$ rather than $\nR^{2d}$, etc.} nothing to do with the operator being bilinear rather than $n$-linear, so it works just as well in our case. There are only two changes that must be made, which are (slightly) past the realm of superficial edits.

To be specific, there were two auxiliary results used, \cite[Proposition B and Lemma C]{dosidis2024boundedness}, which were stated specifically for bilinear operators. We, of course, need their multilinear equivalents. In the case of Lemma C---a bilinear version of Mikhlin's multiplier theorem---the multilinear version is an obvious consequence of \cite[Theorem 5.1]{tomita2010hormander}.\footnote{The process of converting the more general condition on the symbol to the familiar one seen in \cite[Lemma C]{dosidis2024boundedness} is spelled out, e.g., in \cite[Lemma 3.2]{tomita2010hormander}.} Readers interested in such results should check \cite{lee2021hormander}. The multilinear variant of Proposition B is well-known, see \cite[Equation (3.4)]{grafakos2024multilinear}.

\begin{rem}
There is one difference between Theorem \ref{thm:app} and \cite[Theorem 3]{dosidis2024boundedness}. We have had to assume that $1<p_1,\ldots,p_n<\infty$, whereas in \cite{dosidis2024boundedness}, the authors allowed for $1<p_1,p_2\leqslant\infty$. We have given a historical reference for Lemma C, which does not include these edge cases, but the result holds in the more general context, see \cite{lee2021hormander}. On the other hand, Proposition B seems to have not been proven for this edge case in the full multilinear setting. It is likely true in general, however we note that the method used to include the edge cases in Proposition B for the bilinear case does not carry over to the multilinear setting.
\end{rem}

\vspace{+1mm}

\noindent\textbf{Idea of the Proof:} As we will see in Section \ref{sec:vertex}, simply mimicking the argument in \cite{dosidis2024boundedness} ceases to give a sharp result for $n$-linear multipliers with $n > 2$ except in some very special cases. For this reason, we give a novel technique in Section \ref{sec:center}, based on Lemma \ref{lem:changeVariables}, which allows us to remove the shift from a square function, thus permitting us to work solely with shifted maximal functions. This technique gives way for a proof of the sharp result in all cases except those with two or more of $p_1,\ldots,p_n,p'$ being infinity, which corresponds to the situations where the optimal logarithmic exponent \eqref{eq:lambdathm} is $1$. For these remaining points, we apply complex multilinear interpolation, using the special cases from Section \ref{sec:vertex}. For the converse direction, an adjustment of the counterexample given in \cite{dosidis2024boundedness} works, so we use that.

\vspace{+1mm}

\noindent\textbf{Plan of the Paper:} In Section \ref{sec:setup}, we will briefly define the function spaces, with which we will be working, as well as define the shifted operators and state the theorems associated to them that we will need. We end Section \ref{sec:setup} with Lemma \ref{lem:changeVariables}, which will be the crux idea needed later in the paper. Section \ref{sec:proof} begins by reducing the proof of Theorem \ref{thm:main} to the need to establish one estimate, \eqref{eq:logEstimate}, which is the focus of most of the remainder of the paper. In Section \ref{sec:vertex}, we recreate the argument given in \cite{dosidis2024boundedness}, which will give us a sharp result for some special cases, which we interpolate later in Section \ref{sec:interpolation}. Sections \ref{sec:center} and \ref{sec:face} represent the main novelty of the paper; in these, we apply Lemma \ref{lem:changeVariables} along with the shifted maximal function estimates to prove \eqref{eq:logEstimate}, first in the case when $1<p_1,\ldots,p_n,p'<\infty$, and second, when one of the $p$-values is allowed to be infinity. Finally, in Section \ref{sec:converse}, we prove the converse part of Theorem \ref{thm:main} by constructing a counterexample reminiscent of the one in \cite{dosidis2024boundedness}.

\vspace{+1mm}

\noindent\textbf{Notation:} As usual, the conjugate exponent of, say, $p$ will be notated with $p'$, i.e. $\frac{1}{p} + \frac{1}{p'} = 1$. Further, $L_p(\nR^d)$ shall denote the usual space of $p$-integrable functions, $H_1(\nR^d)$ the Hardy space, and $BMO$ the space of functions of bounded mean oscillation. We will also write $L_p(\ell_q)$ for $0<p,q\leqslant\infty$ to mean the space of sequences of complex-valued, Lebesgue measurable functions on $\nR^d$, say $\{g_j\}_j$, with the following (quasi-) norm:
\begin{align*}
    \|\{g_j\}_j\|_{L_p(\ell_q)} &= \left( \int_{\nR^d} \left( \sum_j |g_j(x)|^q \right)^{p/q} \, dx \right)^{1/p} ,
\end{align*}
making the usual modifications if $\max\{p,q\} = \infty$. We will often write $\|g_j\|_{L_p(\ell_q)}$ instead of $\|\{g_j\}_j\|_{L_p(\ell_q)}$. Function spaces in this paper will almost all be on $\nR^d$, so unless otherwise stated, this is assumed and potentially omitted, e.g. we will often write $L_p$ instead of $L_p(\nR^d)$.

Furthermore, we use the notation $A\lesssim B$ to say that there exists $C > 0$ such that $A\leqslant CB$ and $A\lesssim_\gamma B$ to indicate that the aforementioned constant depends on a parameter, $\gamma$. If $A\lesssim B$ and $B\lesssim A$ hold simultaneously, then we write $A\sim B$. We denote the Fourier transform with
\begin{align*}
    \cF f(\xi) = \widehat{f}(\xi) = \int_{\nR^d} f(x) e^{-2\pi i (x,\xi)} \, dx , \quad \xi\in\nR^d ;
\end{align*}
for the inverse, we write $f^\vee(\cdot) = \widehat{f}(-\cdot)$. For the indicator function corresponding to a measurable set, $S$, we write $1_S$.
\vspace{+1mm}

\noindent\textbf{Acknowledgments:} The author would like to express his gratitude towards his Ph.D. advisor, Lenka Slav\'ikov\'a, for her constant support during this project---through all of the wrong turns and insane ideas. Additionally, the author was generously supported during the research for, and writing of, this paper by the Primus research programme PRIMUS/21/SCI/002 of Charles University, the grant no. 23-04720S of the Czech Science Foundation, and the Charles University Research Centre program UNCE/24/SCI/005. 

%% file: Sections/setup.tex
\section{Preliminary Material}\label{sec:setup}

In this section, we will give the necessary ingredients for our proof in Section \ref{sec:proof}. We will be brief, but will give some references with more background and theory development for the interested reader.


\subsection{Shifted Operators}\label{sec:shiftedOperators}

As we have mentioned in the Introduction, the essential tool in our proof will be shifted maximal functions, although we will also use shifted square functions. To give these precisely, we recall the familiar Littlewood-Paley pieces, $\phi,\psi\in\cS(\nR^d)$, satisfying
\begin{align*}
    \supp(\widehat{\phi})\subset \{\xi\in\nR^d : |\xi|\leqslant 2\}, \qquad \widehat{\phi}(\xi) = 1 \quad \text{for} \quad |\xi|\leqslant 1
\end{align*}
and
\begin{align*}
    \supp(\widehat{\psi})\subset \{\xi\in\nR^d :  2^{-1} \leqslant |\xi|\leqslant 2\} , \qquad \sum_{\ell\in\nZ} \widehat{\psi}(2^{-\ell} \xi) = 1
\end{align*}
so long as $\xi\neq 0$ in the latter requirement. These will be fixed throughout the paper. We can define a shifted (and dilated) version of each: for $\ell\in\nZ$ and $y\in\nR^d$,
\begin{align}
    \psi_\ell^y(x) = \psi_\ell(x - 2^{-\ell}y) = 2^{\ell d} \psi(2^\ell x - y) \label{eq:shiftdyadic}
\end{align}
and similarly for $\phi$. Of course, the standard definition of the Hardy space, $H_p$, is that it is those distributions in $\cS'(\nR^d)$, the space of tempered distributions, which, after a suitable refinement, lie in $L_p$; specifically, the $H_p$ norm for $0<p\leqslant \infty$ is given by
\begin{align}
    \|f\|_{H_p(\nR^d)} = \left\| \sup_{\ell\in\nZ} |\phi_\ell * f| \right\|_{L_p(\nR^d)} . \label{eq:hardydefmax}
\end{align}
Equivalently, we have for $0<p<\infty$ that
\begin{align}
    \|f\|_{H_p(\nR^d)} \sim \| \psi_\ell * f \|_{L_p(\ell_2)} . \label{eq:hardydefsquare}
\end{align}
The standard reference for Hardy spaces is \cite[Chapter 3]{stein1993harmonic}, where one can find the definition \eqref{eq:hardydefmax} along with an explanation of the fact that $H_p$ and $L_p$ align when $1<p\leqslant\infty$. The equivalent definition, \eqref{eq:hardydefsquare}, comes from the theory of Triebel-Lizorkin spaces (since $H_p = \dot{F}^0_{p2}$) and it is originally due to J. Peetre \cite{peetre1975spaces}. One may find it in many references, e.g. \cite{frazier1990discrete}.

While \eqref{eq:hardydefsquare} does not hold when $p=\infty$, we can give a definition for $BMO$ by adjusting it. We write $\cD_k$ to denote the dyadic cubes at scale $k\in\nZ$, $\cD$ to be all dyadic cubes, and the quantity $\ell(P)$ to be the side length of $P\in\cD$. With this,
\begin{align*}
    \|f\|_{BMO(\nR^d)} \sim \sup_{P\in\cD} \left( \frac{1}{|P|} \int_P \sum_{\ell = -\log_2 \ell(P)} |(\psi_\ell * f)(x)|^2 \, dx \right)^{1/2} .
\end{align*}
This result also comes from the theory of Triebel-Lizorkin spaces (since $BMO = \dot{F}^0_{\infty 2}$); see \cite{frazier1990discrete}.

\begin{rem}
Obviously this is not the standard definition of $BMO$, but it will do for us; see \cite[Chapter 4]{stein1993harmonic} for background on $BMO$.
\end{rem}

\begin{rem}
Checking \cite{frazier1990discrete,stein1993harmonic}, we see that the above definitions do not depend on the particular choices of $\phi,\psi$.
\end{rem}

By replacing $\phi_\ell$ or $\psi_\ell$ with $\phi_\ell^y$ or $\psi_\ell^y$ (and removing the continuous norms), respectively, we obtain the so-called shifted maximal and square functions: for some $y\in \nR^d$, they are
\begin{align*}
    \sup_{\ell\in\nZ} |(\phi_\ell^y * f)(x)| , \quad \left( \sum_{\ell\in\nZ} |(\psi_\ell^y * f)(x)|^2 \right)^{1/2} ,
\end{align*}
with $f\in\cS'(\nR^d)$. The relationship of these shifted variants to the Lebesgue spaces is elucidated in \cite[Lemma 4]{dosidis2024boundedness} in the case of the square function.

\begin{lem}\label{lem:shiftedsquare}
Suppose $y\in \nR^d$ and $1 \leqslant p < \infty$. Then
\begin{align*}
    \| \{\psi_\ell^y * f\}_{\ell\in\nZ} \|_{L_p(\ell_2)} \lesssim \log(e + |y|)^{\left|\frac12 - \frac{1}{p} \right|} \|f\|_{H_p(\nR^d)} .
\end{align*}
The estimate is sharp in the sense that if $|1/2 - 1/p|$ were smaller, the lemma would not hold. If $p=\infty$, then
\begin{align*}
    \sup_{P\in\cD} \left( \frac{1}{|P|} \int_P \sum_{\ell=-\log_2\ell(P)} |(\psi_\ell^y * f)(x)|^2 \, dx \right)^{1/2} \lesssim \log(e + |y|)^{1/2} \|f\|_{BMO(\nR^d)} .
\end{align*}
Again, the exponent on the logarithm is optimal.
\end{lem}

For the maximal functions, we refer to \cite[Theorem 1.5]{park2024vector}.

\begin{lem}\label{lem:shiftedmax}
Let $1 \leqslant p \leqslant \infty$ and $y\in\nR^d$. Then
\begin{align*}
    \|\{\phi_\ell^y * f\}_{\ell\in\nZ}\|_{L_p(\ell_\infty)} \lesssim \log(e + |y|)^{1/p} \|f\|_{H_p(\nR^d)} .
\end{align*}
Further, $1/p$ as the exponent of the logarithm is sharp in the sense that the estimate would not hold if it were smaller.
\end{lem}

The case $p=\infty$ was not included in \cite[Theorem 1.5]{park2024vector}, but that is because it is obvious and the inequality would just be an equality. Indeed, note the following facts: $\phi_\ell^y * f(x) = \phi_\ell * f(x-2^{-\ell}y)$ and suprema can always be swapped.

\begin{rem}
It was, in fact, a novelty of \cite{dosidis2024boundedness} that the authors proved that the sharp exponent on the logarithm in Lemma \ref{lem:shiftedsquare} is $|\frac12 - \frac{1}{p}|$. The result had previously appeared in the literature, e.g. in \cite[Theorem 5.1]{muscalu2014calderon1} with exponent $1$. After finishing the paper \cite{dosidis2024boundedness}, one of the authors, B. Park, went on to complete the theory of these shifted operators in \cite{park2024vector}. While the optimal logarithm exponent of $1/p$ in Lemma \ref{lem:shiftedmax} was already known prior to \cite{park2024vector} in the literature, see \cite[Chapter 2, Section 5.10]{stein1993harmonic}, it had not been proven before \cite{park2024vector} that $1/p$ is sharp.
\end{rem}


\subsection{The Peetre Maximal Function}

To handle the edge cases in Theorem \ref{thm:main}, we will also need some more technical machinery, the main character of which, will be the Peetre maximal function:
\begin{align*}
    \Mfk_{\sigma,2^k} f(x) = \sup_{z\in\nR^d} \frac{|f(x-z)|}{(1+2^k |z|)^\sigma}
\end{align*}
for $k\in\nZ$ and $\sigma > 0$. This has a long history of application in the area of function spaces and goes back to J. Peetre \cite{peetre1975spaces}.

A nice feature of this maximal function is that there is a Fefferman-Stein-type inequality for it. To state this explicitly, for $r > 0$ we say that $\cE(r)$ is the space of tempered distributions with compactly supported Fourier transforms having support in $\{\xi\in\nR^d : |\xi| \leqslant 2r\}$. Then, assuming $\sigma > d/\min\{p,q\}$, with $0<p<\infty$, $0<q\leqslant\infty$ or $p=q=\infty$, we have
\begin{align}
    \| \{ \Mfk_{\sigma,2^k} f_k \}_{k\in\nZ} \|_{L_p(\ell_q)} \lesssim_{A,p,q} \|\{f_k\}_{k\in\nZ}\|_{L_p(\ell_q)} \label{eq:FeffermanSteinPeetre}
\end{align}
if $f_k\in \cE(A2^k)$ (for some $A>0$, which may be chosen).

Interestingly, one can easily verify that the Peetre maximal function is essentially constant on dyadic cubes in the following sense. If $Q\in\cD_k$, then
\begin{align}
    \sup_{z\in Q} \Mfk_{\sigma,2^k}f(z) \lesssim \inf_{z\in Q} \Mfk_{\sigma,2^k} f(z) . \label{eq:PeetreisConstant}
\end{align}
See \cite{park2021equivalence} (the discussion following Lemma 2.4). Also coming from \cite[Equation (1.3)]{park2021equivalence}, we have the following lemma, which will allow us to incorporate the edge cases into Lemmas \ref{lem:shiftedsquare} and \ref{lem:shiftedmax}.

\begin{lem}\label{lem:auxPeetre}
Suppose $0<p,q\leqslant \infty$, $0<\gamma<1$, $\sigma > d/\min\{p,2\}$, and that $f = \{f_k\}_{k\in\nZ}$ is a sequence with $f_k\in\cE(A 2^k)$ for some $A > 0$. For each $Q\in\cD$, there exists $S_Q\subset Q$, a proper measurable subset, which may depend on $\gamma$ and $f$, such that $|S_Q| > \gamma |Q|$. For $0<p<\infty$ or $p=q=\infty$
\begin{align*}
    \left\| \left\{ \sum_{Q\in \cD_k} \left( \inf_{z\in Q} \Mfk_{\sigma,2^k} f_k(z) \right) 1_{S_Q} \right\}_{k\in\nZ} \right\|_{L_p(\ell_q)} \sim \|f\|_{L_p(\ell_q)}
\end{align*}
and for $p=\infty$,
\begin{align*}
    \left\| \left\{ \sum_{Q\in \cD_k} \left( \inf_{z\in Q} f_k(z) \right) 1_{S_Q} \right\}_{k\in\nZ} \right\|_{L_\infty(\ell_2)}\sim \sup_{P\in\cD} \left( \frac{1}{|P|} \int_P \sum_{k=-\log_2(\ell(P))} |f_k(x)|^2 \, dx \right)^{1/2} .
\end{align*}
\end{lem}

To proceed, we also define
\begin{align*}
    X_p(\nR^d) = \begin{cases}
        H_1(\nR^d) & p=1 , \\
        L_p(\nR^d) & 1<p<\infty , \\
        BMO(\nR^d) & p=\infty .
    \end{cases}
\end{align*}
Combining Lemma \ref{lem:shiftedsquare} with Lemma \ref{lem:auxPeetre}, we have the following, which comes from \cite[Corollary 5]{dosidis2024boundedness}.

\begin{lem}\label{lem:shiftedSquarePeetre}
Let $1\leqslant p \leqslant \infty$, $0<\gamma<1$, $y\in\nR^d$, and $\sigma > d/\min\{p,2\}$. If $f\in X_p(\nR^d)$, there exists $S_Q\subset Q$ for each $Q\in\cD$, a proper, measurable subset (with dependence on $\gamma,y,f$), such that $|S_Q| > \gamma |Q|$ and
\begin{align*}
    \left\| \left\{ \sum_{Q\in\cD_k} \left( \inf_{z\in Q} \Mfk_{\sigma,2^k} (\psi_k^y * f)(z) \right) 1_{S_Q} \right\}_{k\in\nZ} \right\|_{L_p(\ell_2)} \lesssim \log(e + |y|)^{|1/p - 1/2|} \|f\|_{X_p(\nR^d)} .
\end{align*}
\end{lem}

Analogously, we have a statement relating to the shifted maximal function, which, to our knowledge, has not appeared in the literature, although it is an obvious corollary of Lemmas \ref{lem:shiftedmax} and \ref{lem:auxPeetre}.

\begin{lem}\label{lem:shiftedMaxPeetre}
Let $1\leqslant p \leqslant \infty$, $0<\gamma<1$, $y\in\nR^d$, and $\sigma > d/\min\{p,2\}$. Then for $f\in H_p(\nR^d)$ and each $Q\in\cD$, there exists $S_Q\subset Q$, a proper measurable subset (depending on $\gamma,y,f$), satisfying the condition $|S_Q| > \gamma |Q|$, such that
\begin{align*}
    \left\| \left\{ \sum_{Q\in\cD_k} \left( \inf_{z\in Q} \Mfk_{\sigma,2^k} (\phi_k^y * f)(z) \right) 1_{S_Q} \right\}_{k\in\nZ} \right\|_{L_p(\ell_\infty)} \lesssim \log(e + |y|)^{1/p} \|f\|_{H_p(\nR^d)} .
\end{align*}
\end{lem}


\subsection{A Trick with Change of Variables}

The final tool that we require is a simple lemma, which says we can remove a shift from one function, placing it on others, if we are in an $L_2(\ell_2)$-norm.

\begin{lem}\label{lem:changeVariables}
Suppose $g_1,\ldots,g_m$ are measurable. Define their shifted dyadic dilates as in \eqref{eq:shiftdyadic}: for $k=1,\ldots,m$,
\begin{align*}
    g_{k\ell}^{y_k} = g_k(x - 2^{-\ell} y_k) = 2^{\ell d} g_k(2^\ell x - y_k)
\end{align*}
with $y_k\in \nR^d$. Then for any $k_0 \in \{1,\ldots,m\}$,
\begin{align*}
    \left\| \left\{ \prod_{k=1}^m g_{k\ell}^{y_k} \right\}_{\ell\in\nZ} \right\|_{L_2(\ell_2)} = \left\| \left\{ g_{k_0 \ell} \prod_{\substack{k=1 \\ k\neq k_0}}^m g_{k\ell}^{y_k-y_{k_0}} \right\}_{\ell\in\nZ} \right\|_{L_2(\ell_2)} .
\end{align*}
In words, we can remove the shift from any one of the $g_k$ and place it onto the others.
\end{lem}
\begin{proof}
This is nothing more than a change of variables. Indeed, because we have an $L_2(\ell_2)$ norm the square and square root cancel, allowing us to swap the sum and integral:
\begin{align*}
    \left\| \left\{ \prod_{k=1}^m g_{k\ell}^{y_k} \right\}_{\ell\in\nZ} \right\|_{L_2(\ell_2)}^2 = \int_{\nR^d} \sum_{\ell\in\nZ} \left| \prod_{k=1}^m g_{k\ell}^{y_k}(x) \right|^2 \, dx = \sum_{\ell\in\nZ} \int_{\nR^d} \left| \prod_{k=1}^m g_{k\ell}^{y_k}(x) \right|^2 \, dx .
\end{align*}
Perform the change of variables $x - 2^{-\ell}y_{k_0} \mapsto x$ to get
\begin{align*}
    \sum_{\ell\in\nZ} \int_{\nR^d} \left| g_{k_0\ell}(x) \prod_{\substack{k=1 \\ k\neq k_0}}^m g_{k\ell}^{y_k - y_{k_0}}(x) \right|^2 \, dx .
\end{align*}
Finally, we can return to the $L_2(\ell_2)$ norm by reversing the steps at the beginning of the proof.
\end{proof}

\begin{rem}
The reader may be wondering what is special about $L_2(\ell_2)$ and the answer is that there is nothing special about it. The salient detail was that the square and square root canceled, allowing us to swap the sum and integral. This would have worked just as well with powers of $p$ and $1/p$, i.e. if we had used $L_p(\ell_p)$ for $0<p<\infty$.
\end{rem}

%% file: Sections/proof.tex
\section{Proof of the Positive Result}\label{sec:proof}

This section will be divided into five parts. In the first, we will reduce the proof of Theorem \ref{thm:main} to the verification of one estimate, \eqref{eq:logEstimate}. This part simply replicates what was done in \cite{dosidis2024boundedness}, making edits to fit the $n$-linear---rather than bilinear---setting. The remaining four parts distinguish four cases, where the validity of \eqref{eq:logEstimate} will be shown. The first will imitate the method of \cite{dosidis2024boundedness}, which now only gives the sharp result in very particular cases. The following two will follow a new method using shifted maximal functions, which is the main novelty of this paper. In the final part, we apply complex interpolation to clean up the remaining contingencies.


\subsection{Reducing the Proof to One Estimate}\label{sec:proofSetup}

Define the $(n+1)$-linear form
\begin{align}
    \Lambda(f_1,\ldots,f_n,f_{n+1}) = \int_{\nR^d} T(f_1,\ldots,f_n)(x) f_{n+1}(x) \, dx \label{eq:Lambda}
\end{align}
for $f_1,\ldots,f_{n+1}\in\cS_0(\nR^d)$. Given that our goal is \eqref{eq:maineq}, we observe that it suffices to demonstrate
\begin{align}
    |\Lambda(f_1,\ldots,f_n,f_{n+1})| \lesssim_{p_1,\ldots,p_n,d,n} D_\lambda(K) \|f_1\|_{H_{p_1}(\nR^d)} \cdots \|f_n\|_{H_{p_n}(\nR^d)} \|f_{n+1}\|_{H_{p'}(\nR^d)} . \label{eq:goal}
\end{align}
Before continuing, however, we define the transpose operators $T^1,\ldots,T^n$ in the usual way: the $j$-th transpose operator is given by swapping $f_j$ with $f_{n+1}$ in \eqref{eq:Lambda}:
\begin{align*}
    \Lambda(f_1,\ldots,f_n,f_{n+1}) = \int_{\nR^d} T^j(f_1,\ldots,f_{j-1},f_{n+1},f_{j+1},\ldots,f_n)(x) f_j(x) \, dx
\end{align*}
for $1\leqslant j\leqslant n$ with the obvious modifications to the above expression if $j=1$ or $j=n$. Through a simple change of variables, the kernels of each $T^j$, denoted $K^j$ for $1\leqslant j\leqslant n$, can be related back to $K$ with
\begin{align}
    K^j(y_1,\ldots,y_n) &= K(y_1-y_j,y_2-y_j,\ldots,y_{j-1}-y_j,-y_j,y_{j+1}-y_j,\ldots,y_n-y_j) , \label{eq:transposeKernel}
\end{align}
again, with the obvious modification if $j=1$ or $j=n$. These are really just rotations of $K$, so
\begin{align}
    D_\lambda(K) \sim_d D_\lambda(K^j) \label{eq:simDlambda}
\end{align}
for each $j$. Additionally, their Fourier transforms are all still supported in $|(\xi_1,\ldots,\xi_n)|\sim_d 1$---the unit annulus---just like that of $K$. 

We now proceed with the proof. It will be convenient to write $\xi = (\xi_1,\ldots,\xi_n)\in\nR^{nd}$ and $y = (y_1,\ldots,y_n)\in\nR^{nd}$. Use Fourier inversion to rewrite
\begin{align}
    \Lambda(f_1,\ldots,f_{n+1}) &= \sum_{\ell\in\nZ} \int_{\nR^d} \cF(K_\ell * (f_1\otimes \cdots \otimes f_n))^\vee (x,\ldots,x) f_{n+1}(x) \, dx \nonumber \\
    &= \sum_{\ell\in\nZ} \int_{\nR^{nd}} \widehat{K}(2^{-\ell} \xi_1, \ldots, 2^{-\ell}\xi_n) \left( \prod_{k=1}^n \widehat{f_k}(\xi_k) \right) \widehat{f}_{n+1}(-\xi_1-\cdots -\xi_n) \, d\xi , \label{eq:fourierinv}
\end{align}
since
\begin{align*}
    \int_{\nR^d} f_{n+1}(x) e^{2\pi i ((x,\ldots,x),(\xi_1,\ldots,\xi_n))} \, dx = \int_{\nR^d} f_{n+1}(x) e^{-2\pi i(x,-\xi_1-\cdots -\xi_n)} \, dx = \widehat{f}_{n+1}(-\xi_1-\cdots -\xi_n) .
\end{align*}

With this expression for $\Lambda$, we would like to take advantage of the support of $\widehat{K}$. Indeed, with it vanishing outside of the unit annulus, it must be that at least two of $\xi_1,\ldots,\xi_n$, and $-\xi_1-\cdots -\xi_n$ are bounded away from zero in its support. In accordance with that, we recall the two Schwartz functions, $\phi,\psi\in\cS(\nR^n)$ from Section \ref{sec:setup}. Using these and an argument with a smooth partition of unity, we can decompose $\widehat{K}$ into $\binom{n+1}{2}$ pieces, each of the form
\begin{align}
    \widehat{K}(\xi_1,\ldots,\xi_n) \left( \prod_{k=1}^n \widehat{\Phi}_k(\xi_k) \right) \widehat{\Phi}_{n+1}(-\xi_1-\cdots -\xi_n) , \label{eq:Kdecomp}
\end{align}
where two of the $\Phi_j$ are $\psi$ and the others are $\phi$.

Each part of $\widehat{K}$ in our decomposition gives rise to an $(n+1)$-linear form reminiscent of $\Lambda$. To progress, we need to be able to talk about pairs $(s,t)\in\nN^2$ without ambiguity. We let $\cI$ be the set of all pairs $(s,t)$ with $1\leqslant s < t \leqslant n+1$. For any $(s,t)\in\cI$, the operator $\Lambda_{(s,t)}$ is, then, defined as
\begin{align}
    &\Lambda_{(s,t)}(f_1,\ldots,f_{n+1}) \notag \\
    &\hspace{+5mm}= \sum_{\ell\in\nZ} \int_{\nR^{nd}} \widehat{K}(2^{-\ell}\xi) \left(\prod_{k=1}^n \widehat{\Phi}_{k\ell}(\xi_k) \widehat{f_k}(\xi_k)\right) \widehat{\Phi}_{n+1,\ell}(-\xi_1-\cdots-\xi_n) \widehat{f}_{n+1}(-\xi_1-\cdots-\xi_n) \, d\xi , \label{eq:changeTranspose}
\end{align}
where $\Phi_k = \psi$ for $k = s,t$ and $\Phi_k = \phi$ otherwise. We will be using the notation $\Phi_{k\ell}^{y_k}(x) = 2^{\ell d} \Phi_k(2^\ell x - y_k)$ and $\Phi_{k\ell}^0 = \Phi_{k\ell}$. Based on our progress thus far,
\begin{align*}
    |\Lambda(f_1,\ldots,f_{n+1})| \lesssim \sum_{(s,t)\in\cI} |\Lambda_{(s,t)} (f_1,\ldots,f_{n+1})| .
\end{align*}

Returning to our transpose operators, we can rearrange at will which of the $\Phi_j$ is evaluated at $-\xi_1-\cdots -\xi_n$: apply \eqref{eq:transposeKernel} to see
\begin{align*}
    \widehat{K}(\xi_1,\ldots,\xi_n) = \widehat{K^j}(\xi_1,\ldots,\xi_{j-1},-\xi_1-\cdots -\xi_n,\xi_{j+1},\ldots ,\xi_n)
\end{align*}
for $1\leqslant j\leqslant n$ (with the obvious change if $j=1$ or $j=n$). Using this in \eqref{eq:changeTranspose}, we can then do a simple change of variables to move the argument $-\xi_1-\cdots -\xi_n$ around. In words, we can choose at will which of the functions in \eqref{eq:changeTranspose} is evaluated at $-\xi_1-\cdots-\xi_n$; in $\Lambda_{(s,t)}$, this point will be, say, $f_\tau$ for some $1\leqslant \tau \leqslant n+1$ (we will choose $\tau$ based on the context). For simplicity of notation, we write $p_{n+1} = p'$ and $K^{n+1} = K$. Therefore,
\begin{align*}
    &\Lambda_{(s,t)}(f_1,\ldots,f_{n+1}) \\
    &\hspace{+5mm}= \sum_{\ell\in\nZ} \int_{\nR^{nd}} \widehat{K^\tau}(2^{-\ell}\xi) \left( \prod_{\substack{1\leqslant k\leqslant n+1 \\ k\neq\tau}} \widehat{\Phi}_{k\ell}(\xi_k) \widehat{f_k}(\xi_k) \right) \widehat{\Phi}_{\tau\ell}(-\xi_1-\cdots-\xi_n) \widehat{f}_\tau(-\xi_1-\cdots-\xi_n) \, d\xi .
\end{align*}
Applying Fourier inversion again, i.e. proceeding as in \eqref{eq:fourierinv}, but in reverse, we find
\begin{align*}
    &\Lambda_{(s,t)}(f_1,\ldots,f_{n+1}) = \sum_{\ell\in\nZ} \int_{\nR^{nd}} \int_{\nR^d} K^\tau (y) \left(\prod_{\substack{1\leqslant k\leqslant n+1 \\ k\neq\tau}} (\Phi_{k\ell}^{y_k} * f_k)(x) \right) (\Phi_{\tau\ell} * f_\tau)(x) \, dx dy .
\end{align*}
Above, if $\tau\neq n+1$, we have written $y_{n+1} := y_\tau$.

\begin{rem}
It is at this point above that we see the inevitability of shifted square and maximal functions, given our approach. Indeed, if we take a convolution-type operator and decompose it in Fourier space while multiplying in Littlewood-Paley pieces, then we have introduced a second convolution in physical space, which is reflected with the addition of the shifts on each $\Phi_k$.
\end{rem}

From here, we proceed by first bounding the above without the integral over $y$ (or the kernel). If we can show
\begin{align}
    &\left| \sum_{\ell\in\nZ} \int_{\nR^d} \left(\prod_{\substack{1\leqslant k\leqslant n+1 \\ k\neq\tau}} (\Phi_{k\ell}^{y_k} * f_k)(x) \right) (\Phi_{\tau\ell} * f_\tau)(x) \, dx \right| \lesssim \log(e + |y|)^{\lambda_{(s,t)}} \prod_{k=1}^{n+1} \|f_k\|_{H_{p_k}} , \label{eq:logEstimate}
\end{align}
where
\begin{align}
    \lambda_{(s,t)} = \sum_{k\in\{1,\ldots,n+1\}\smallsetminus \{s,t\}} \frac{1}{p_k} \label{eq:lambdast} ,
\end{align}
then we can bound each $\lambda_{(s,t)}$ above by the right hand side of \eqref{eq:lambdathm}, call this $\lambda$. Reintroducing the integral over $y$ and using the definition of $D_\lambda(K^\tau)$,
\begin{align*}
    |\Lambda_{(s,t)}(f_1,\ldots,f_{n+1})| \lesssim D_{\lambda}(K^{\tau}) \prod_{k=1}^{n+1} \|f_k\|_{H_{p_k}} .
\end{align*}
Therefore, combining each estimate (also using \eqref{eq:simDlambda}), we see that
\begin{align*}
    |\Lambda(f_1,\ldots,f_{n+1})| \lesssim \sum_{(s,t)\in\cI} |\Lambda_{(s,t)}(f_1,\ldots,f_{n+1})| \lesssim D_{\lambda}(K) \prod_{k=1}^{n+1} \|f_k\|_{H_{p_k}} ,
\end{align*}
as desired. It, therefore, remains to establish \eqref{eq:logEstimate}.


\subsection{The Vertex Points}\label{sec:vertex}

We imagine this as proving the result by giving it for each point of the form $(\frac{1}{p_1},\ldots,\frac{1}{p_n},\frac{1}{p'})$; of course, we are assuming that the coordinates of this point sum to one. Such points make a shape\footnote{It is the lower half of a cube---below the diagonal.} in $\nR^n$ (really, $\nR^{n+1}$, but the final point is always determined by the rest) and the goal of this section is to prove the theorem for its vertices, i.e. $(1,0,\ldots,0), (0,1,0,\ldots,0), \ldots , (0,0,\ldots,0,1)$. It is, actually, no extra work to just prove a result for an arbitrary point $(\frac{1}{p_1},\ldots,\frac{1}{p_n},\frac{1}{p'})$, which is only sharp at the vertices.

Recall that, in \eqref{eq:logEstimate}, the pair $(s,t)$ corresponds to those $\Phi_k$ which are equal to $\psi$ and that $\tau$ is the distinguished component, which has no shift attached to it. Assume that $\tau \neq s,t$ (we will discuss what happens when $\tau$ is $s$ or $t$ in Remark \ref{rem:tau=st}). The particular choice of $\tau$, we leave for later.

First, we set up our use of Lemmas \ref{lem:shiftedSquarePeetre} and \ref{lem:shiftedMaxPeetre}. Take $\sigma > d$ and define for $z\in \nR^d$
\begin{align}
    \Qfk_{\sigma,\ell}^z f = \Mfk_{\sigma,2^\ell}(\psi_\ell^z * f) , \quad \Rfk_{\sigma,\ell}^z f = \Mfk_{\sigma,2^\ell}(\phi_\ell^z * f) . \label{eq:QRfk}
\end{align}
We will, additionally, write $\Rfk_{\sigma,\ell} f = \Rfk_{\sigma,\ell}^0 f$ and $\Qfk_{\sigma,\ell} f = \Qfk_{\sigma,\ell}^0 f$. Applying Lemmas \ref{lem:shiftedSquarePeetre} and \ref{lem:shiftedMaxPeetre} to each function, $f_k$ with $k\neq \tau$, we pick the sequences $(S_Q^k)_{Q\in\cD}$ of measurable sets $S_Q^k \subset Q$, where the inclusion is proper, such that
\begin{align*}
    |S_Q^k| > \left( 1 - \frac{1}{2n} \right) |Q| 
\end{align*}
for each $k\neq \tau$. The condition above guarantees
\begin{align}
    \left| \bigcap_{k\neq\tau} S_Q^k \right| \geqslant \frac12 |Q| \label{eq:intersec}
\end{align}
by a simple iterative calculation. Call the intersection $S_Q$, i.e. $S_Q = \bigcap S_Q^k$. Further, we have
\begin{align}
    \left\| \left\{ \sum_{Q\in\cD_\ell} \left( \inf_{z\in Q} \Qfk_{\sigma,\ell}^{y_k} f_k(z) \right) 1_{S_Q^k} \right\}_{\ell\in\nZ} \right\|_{L_{p_k}(\ell_2)} \lesssim \log(e+|y_k|)^{|1/p_k - 1/2|} \|f_k\|_{H_{p_k}(\nR^d)} \label{eq:squareInProof}
\end{align}
for $k=s,t$\footnote{Looking at Lemma \ref{lem:shiftedSquarePeetre}, if $p_s$ or $p_t$ is $\infty$, then we would get the $BMO$ norm on the right hand side, but this is, of course, bounded by the $H_\infty = L_\infty$ norm.} and
\begin{align}
    \left\| \left\{ \sum_{Q\in\cD_\ell} \left( \inf_{z\in Q} \Rfk_{\sigma,\ell}^{y_k} f_k(z) \right) 1_{S_Q^k} \right\}_{\ell\in\nZ} \right\|_{L_{p_k}(\ell_\infty)} \lesssim \log(e+|y_k|)^{1/p_k} \|f_k\|_{H_{p_k}(\nR^d)} \label{eq:maxInProof}
\end{align}
for $k\in \{1,\ldots,n+1\}\smallsetminus \{s,t,\tau\}$ (if $n=2$, this is void).

The left hand side of \eqref{eq:logEstimate} is bounded by
\begin{align*}
    &\sum_{\ell\in\nZ} \sum_{Q\in \cD_\ell} \int_Q \Qfk_{\sigma,\ell}^{y_s} f_s(x) \Qfk_{\sigma,\ell}^{y_t} f_t(x) \left( \prod_{k\neq s,t,\tau} \Rfk_{\sigma,\ell}^{y_k} f_k(x) \right) \Rfk_{\sigma,\ell}f_\tau(x) \, dx \\
    &\hspace{+3mm} \lesssim \sum_{\ell\in\nZ} \sum_{Q\in\cD_\ell} \left( \inf_{z\in Q} \Qfk_{\sigma,\ell}^{y_s} f_s(z) \right) \left( \inf_{z\in Q} \Qfk_{\sigma,\ell}^{y_t}f_t(z) \right) \left( \prod_{k\neq s,t,\tau} \inf_{z\in Q} \Rfk_{\sigma,\ell}^{y_k} f_k(z) \right) \left( \inf_{z\in Q} \Rfk_{\sigma,\ell} f_\tau(z) \right) |Q| .
\end{align*}
We have used \eqref{eq:PeetreisConstant} in the second estimate above. Appealing to \eqref{eq:intersec}, this is dominated by
\begin{align*}
    \sum_{\ell\in\nZ} \int_{\nR^d} \Rfk_{\sigma,\ell} f_\tau (x) \sum_{Q\in \cD_\ell} \left[ \left( \inf_{z\in Q} \Qfk_{\sigma,\ell}^{y_s} f_s(z) \right) \left( \inf_{z\in Q} \Qfk_{\sigma,\ell}^{y_t}f_t(z) \right) \left( \prod_{k\neq s,t,\tau} \inf_{z\in Q} \Rfk_{\sigma,\ell}^{y_k} f_k(z) \right) \right] 1_{S_Q}(x) \, dx .
\end{align*}
Applying H\"older's inequality then gives the bounds\footnote{Recall that we have let $p_{n+1} = p'$.}
\begin{align*}
    &\int_{\nR^d} \|\{\Rfk_{\sigma,\ell}f_\tau (x)\}_{\ell\in\nZ}\|_{\ell_\infty} \prod_{k\neq s,t,\tau}  \left\| \left\{ \sum_{Q\in \cD_\ell} \left( \inf_{z\in Q} \Rfk_{\sigma,\ell}^{y_k}f_k(z) \right) 1_{S_Q^k}(x) \right\}_{\ell\in\nZ} \right\|_{\ell_\infty} \\
    &\hspace{+3cm}\cdot \prod_{k=s,t}  \left\| \left\{ \sum_{Q\in \cD_\ell} \left( \inf_{z\in Q} \Qfk_{\sigma,\ell}^{y_k}f_k(z) \right) 1_{S_Q^k}(x) \right\}_{\ell\in\nZ} \right\|_{\ell_2} \, dx \\
    &\hspace{+3mm}\leqslant \|\{\Rfk_{\sigma,\ell}f_\tau\}_{\ell\in\nZ}\|_{L_{p_\tau}(\ell_\infty)} \prod_{k\neq s,t,\tau} \left\| \left\{ \sum_{Q\in \cD_\ell} \left( \inf_{z\in Q} \Rfk_{\sigma,\ell}^{y_k}f_k(z) \right) 1_{S_Q^k} \right\}_{\ell\in\nZ} \right\|_{L_{p_k}(\ell_\infty)} \\
    &\hspace{+3cm}\cdot \prod_{k=s,t} \left\| \left\{ \sum_{Q\in \cD_\ell} \left( \inf_{z\in Q} \Qfk_{\sigma,\ell}^{y_k}f_k(z) \right) 1_{S_Q^k} \right\}_{\ell\in\nZ} \right\|_{L_{p_k}(\ell_2)} .
\end{align*}
Finally, using \eqref{eq:squareInProof}, \eqref{eq:maxInProof} for the two products and \eqref{eq:FeffermanSteinPeetre} followed by the definition of $H_p$ for the $f_\tau$ term, we see that this final display is bounded by the right hand side of \eqref{eq:logEstimate}, except with $\lambda_{(s,t)}'$, where
\begin{align}
    \lambda_{(s,t)}' &= \left|\frac12 - \frac{1}{p_s} \right| + \left|\frac12 - \frac{1}{p_t} \right| + \sum_{k\neq s,t,\tau} \frac{1}{p_k} \label{eq:lambdst'}
\end{align}
instead of $\lambda_{(s,t)}$ from \eqref{eq:lambdast}. By a simple computation, this is bounded by
\begin{align}
    \lambda' = \max\left\{ \frac{1}{p_1} , \ldots , \frac{1}{p_n} , \frac{1}{p'} \right\} + 2 \left( \sum_{k=2}^{n-1} \mx{k} \left\{ \frac{1}{p_1} , \ldots , \frac{1}{p_n} , \frac{1}{p'} \right\} \right) \label{eq:lambda'}
\end{align}
for all pairs $(s,t)\in\cI$. Indeed, if $\frac{1}{p_k} < \frac12$ for each $1\leqslant k\leqslant n+1$, then \eqref{eq:lambda'} is clearly the largest that \eqref{eq:lambdst'} can be. The maximum, in this case, is achieved when $\frac{1}{p_s}$ and $\frac{1}{p_t}$ are the two smallest coordinates in $(\frac{1}{p_1},\ldots,\frac{1}{p_n},\frac{1}{p'})$ and $\tau$ is chosen so that $\frac{1}{p_\tau}$ is as large as possible (we want to exclude the biggest one), i.e. the largest of the $\frac{1}{p_k}$ with $k\neq s,t$. If $\frac{1}{p_u} = \frac{1}{p_v} = \frac12$ for some $1\leqslant u,v\leqslant n+1$, then the claim is obvious. If, on the other hand, we have $\frac{1}{p_u}\geqslant \frac12$ for some $1\leqslant u\leqslant n+1$ and $\frac{1}{p_k} < \frac12$ for $k\neq u$, then there are two possibilities. If $u$ is neither $s$ nor $t$, then we proceed as before by choosing $\tau = u$. If $u = s$ (or, similarly, if $u = t$), then this is not the pair $(s,t)$ that maximizes $\lambda_{(s,t)}'$ over all $(s,t)\in\cI$. Observe, in this case,
\begin{align*}
    \lambda_{(u,t)}' = \left( \frac{1}{p_u} - \frac12 \right) + \left( \frac12 - \frac{1}{p_t} \right) + \sum_{k\neq u,t,\tau} \frac{1}{p_k} = \left( \sum_{k\neq t,\tau} \frac{1}{p_k} \right) - \frac{1}{p_t} .
\end{align*}
The largest this can be is when $\frac{1}{p_t}$ is as small as possible and $\frac{1}{p_\tau}$ is the second smallest. Looking at this sub-case, the final expression above is equal to
\begin{align*}
    \max\left\{ \frac{1}{p_1} , \ldots , \frac{1}{p_n} , \frac{1}{p'} \right\} + 1 \left( \sum_{k=2}^{n-1} \mx{k} \left\{ \frac{1}{p_1} , \ldots , \frac{1}{p_n} , \frac{1}{p'} \right\} \right) - \min\left\{ \frac{1}{p_1} , \ldots , \frac{1}{p_n} , \frac{1}{p'} \right\} ,
\end{align*}
which is smaller than $\lambda'$; in fact, it is always smaller than the optimal $\lambda$ in \eqref{eq:lambdathm}.

It only, now, remains to note that $\lambda'$ above equals the right hand side of \eqref{eq:lambdathm} if $(\frac{1}{p_1} , \ldots , \frac{1}{p_n} , \frac{1}{p'})$ is a vertex point, which means we have proven the optimal result for these vertices.

\begin{rem}
The above argument is a na\"ive reproduction of the proof in \cite{dosidis2024boundedness} of the bilinear case; it is interesting to see that in the general $n$-linear case, with $n > 2$, it is only sharp at the vertex points.
\end{rem}

\begin{rem}\label{rem:tau=st}
One could also redo the above proof choosing $\tau$ to be $s$ or $t$, say $\tau = t$. Then we would end with
\begin{align*}
    \lambda_{(s,t)}'' = \left| \frac12 - \frac{1}{p_s} \right| + \sum_{k\neq s,t} \frac{1}{p_k} .
\end{align*}
This does turn out to be useful; in fact, we will use it at the end of Section \ref{sec:center}.
\end{rem}


\subsection{The Non-Edge Cases}\label{sec:center}

With the vertices done, we turn our attention to points of the form $(\frac{1}{p_1},\ldots,\frac{1}{p_n},\frac{1}{p'})$, where $1<p_1,\ldots,p_n,p'<\infty$. The argument in this section, when combined with the opening steps of the previous section, essentially works when one of the $1/p_k$ values is zero, but will break down when two or more are. We will describe these small adjustments in the next section.

On the left hand side of \eqref{eq:logEstimate}, we will, again, take $(s,t)\in \cI$ arbitrary; this time, however, we assume that $\tau = s$ or $\tau = t$. For the purposes of this section, it does not matter if we choose $\tau = s$ or $\tau = t$, so we will let $\tau = t$ for concreteness. In words, this means that one of the square functions will have no shift attached to it; with the other, our goal is to move the shift off of it.

The following argument has the unfortunate drawback of, in spirit, being simple, but in rigorous presentation, being somewhat complex. We describe the spirit first to orient ourselves. As mentioned, we want to move the shift off of the square function and onto some of the maximal functions; for this we will apply Lemma \ref{lem:changeVariables}. Now, in order to apply this lemma, we have to create an $L_2(\ell_2)$ norm, which will require us to group the $1/p_k$ values together in such a way that they add to exactly $1/2$; to do so, we will also have to split one of these values in two. Once we have created the $L_2(\ell_2)$ norm through H\"older's inequality, it will only remain to apply the lemma and proceed.

We distinguish two cases: in the first, it holds that $\frac{1}{p_k} < \frac12$ for $k=1,\ldots,n+1$.\footnote{Recall that we write $p' = p_{n+1}$.} Obviously, the other case will be the negation of this.

In the former case, since each reciprocal is less than $\frac12$ and we wish to remove the shift from $\psi^{y_s}*f_s$, we define the subset\footnote{Note that $1/p_t$ is forbidden in the set below, because the $L_{p_t}$ norm will also come from a square function. Moving the shift from one square function to the other would defeat the purpose of the process below. We only exclude $1/p_s$ to distinguish it from the others.}
\begin{align*}
    J_0 \subset \left\{ 1,\ldots,n+1 \right\}\smallsetminus \{s,t\} ,
\end{align*}
such that $\sum_{k\in J_0\cup \{s\}} \frac{1}{p_k} < \frac12$, and there exists $\alpha\in (J_0)^c$ such that $\sum_{k\in J_0\cup \{\alpha,s\}} \frac{1}{p_k} > \frac12$.\footnote{We could explicitly construct $J_0$ by reading from left to right. For example, if $s = n$, then check if $1/p_n + 1/p_1 > 1/2$. If it is, then $J_0$ is empty and $\alpha = 1$. If not, then check if $1/p_n + 1/p_1 + 1/p_2 > 1/2$. If it is, then $J_0 = \{1\}$ and $\alpha = 2$, etc. Obviously this process will terminate since $1/p_t < 1/2$, so the sum of the rest is larger than $1/2$.} It is possible that there exists a subset of $\{1/p_1,\ldots,1/p_{n+1}\}\smallsetminus \{1/p_s,1/p_t\}$ that, when combined with $1/p_s$, has sum exactly equal to $\frac12$. The steps below will, then, just become simpler: we can either choose $J_0$ to be such that $\sum_{k\in J_0\cup \{s\}} \frac{1}{p_k} = \frac12$ and avoid choosing $\alpha$ altogether or do the same as below and just let $\alpha\in (J_0)^c$ be arbitrary (and, below, $\gamma = 0$). The edits to the argument to come, in this case, are obvious.

What we do now is break
\begin{align*}
    \frac{1}{p_\alpha} = \frac{1}{q_0} + \frac{1}{q_1} ,
\end{align*}
where this decomposition is chosen so that $1<q_0,q_1<\infty$ and
\begin{align}
    \left( \sum_{k\in J_0\cup \{s\}} \frac{1}{p_k} \right) + \frac{1}{q_0} = \frac12 . \label{eq:exactly1/2}
\end{align}
For convenience, set $J_1 = (J_0)^c \smallsetminus \{\alpha\}$.\footnote{This also does not include $s$ or $t$!} Returning with this setup to our goal, the left hand side of \eqref{eq:logEstimate} is bounded by
\begin{align*}
    &\sum_{\ell\in\nZ} \int_{\nR^d} \left( \prod_{k\in J_1} |(\phi_\ell^{y_k} * f_k)(x)| \right) |(\phi_\ell^{y_\alpha} * f_\alpha)(x)|^{1-\gamma} \\
    &\hspace{+15mm}\cdot\left| (\psi_\ell^{y_s} * f_s)(x) |(\phi_\ell^{y_\alpha} * f_\alpha)|^\gamma \prod_{k\in J_0} (\phi_\ell^{y_k} * f_k)(x) \right| |(\psi_\ell * f_t)(x)| \, dx ,
\end{align*}
where $0<\gamma<1$ is chosen so that $q_0\gamma = p_\alpha$ (and, therefore, $q_1 (1-\gamma) = p_\alpha$). H\"older's inequality shows this is bounded by
\begin{align*}
    &\int_{\nR^d} \left( \prod_{k\in J_1} \sup_{\ell\in\nZ} |(\phi_\ell^{y_k} * f_k)(x)| \right) \sup_{\ell\in\nZ} |(\phi_\ell^{y_\alpha} * f_\alpha)(x)|^{1-\gamma} \\
    &\hspace{+15mm}\cdot\left( \sum_{\ell\in\nZ} \left| (\psi_\ell^{y_s} * f_s)(x) |(\phi_\ell^{y_\alpha} * f_\alpha)|^\gamma \prod_{k\in J_0} (\phi_\ell^{y_k} * f_k)(x) \right|^2 \right)^{\frac12} \left( \sum_{\ell\in\nZ} |(\psi_{\ell} * f_t)(x)|^2 \right)^{\frac12} \, dx \\
    &\hspace{+5mm}\leqslant \left( \prod_{k\in J_1} \|\phi_\ell^{y_k} * f_k\|_{L_{p_k}(\ell_\infty)} \right) \|\phi_\ell^{y_\alpha} * f_\alpha\|_{L_{p_\alpha}(\ell_\infty)}^{1-\gamma} \\
    &\hspace{+15mm}\cdot\left\| (\psi_\ell^{y_s} * f_s)(x) |(\phi_\ell^{y_\alpha} * f_\alpha)|^\gamma \prod_{k\in J_0} (\phi_\ell^{y_k} * f_k)(x) \right\|_{L_2(\ell_2)} \|\psi_\ell * f_t\|_{L_{p_t}(\ell_2)} .
\end{align*}
Above, we have used that
\begin{align*}
    \||\phi_\ell^{y_\alpha} * f_\alpha|^{1-\gamma}\|_{L_{q_1}(\ell_\infty)} = \|\phi_\ell^{y_\alpha} * f_\alpha\|_{L_{p_\alpha}(\ell_\infty)}^{1-\gamma} ;
\end{align*}
this just follows from the definition of $L_p$ norms and the fact that $(1-\gamma) q_1 = p_\alpha$.\footnote{We also had to swap a supremum and $\gamma$ power function; this is obviously possible, because the $\gamma$ power function is increasing.}

Having created the $L_2(\ell_2)$ norm, we apply Lemma \ref{lem:changeVariables} to see that this term in our expression above is equal to
\begin{align*}
    \left\| (\psi_\ell * f_s)(x) |(\phi_\ell^{y_\alpha - y_s} * f_\alpha)|^\gamma \prod_{k\in J_0} (\phi_\ell^{y_k - y_s} * f_k)(x) \right\|_{L_2(\ell_2)} .
\end{align*}
Therefore, using H\"older's inequality again, the whole expression is bounded by
\begin{align*}
    &\left( \prod_{k\in J_1} \|\phi_\ell^{y_k} * f_k\|_{L_{p_k}(\ell_\infty)} \right) \|\phi_\ell^{y_\alpha} * f_\alpha\|_{L_{p_\alpha}(\ell_\infty)}^{1-\gamma} \|\psi_\ell * f_s\|_{L_{p_s}(\ell_2)} \|\phi_\ell^{y_\alpha - y_s} * f_\alpha\|_{L_{p_\alpha}(\ell_\infty)}^\gamma \\
    &\hspace{+15mm}\cdot\left( \prod_{k\in J_0} \|\phi_\ell^{y_k - y_s} * f_k\|_{L_{p_k}(\ell_\infty)} \right) \|\psi_\ell * f_t\|_{L_{p_t}(\ell_2)} .
\end{align*}
Finally using Lemma \ref{lem:shiftedmax} and the classical result for square functions (Lemma \ref{lem:shiftedsquare} with $y=0$ or \eqref{eq:hardydefsquare}), this is bounded, up to a constant, by
\begin{align}
    &\left( \prod_{k\neq s,t,\alpha} \log(e + |y|)^{\frac{1}{p_k}} \|f_k\|_{L_{p_k}} \right) \left( \log(e + |y|)^{\frac{1-\gamma}{p_\alpha}} \|f_\alpha\|_{L_{p_\alpha}}^{1-\gamma} \log(e+|y|)^{\frac{\gamma}{p_\alpha}} \|f_\alpha\|_{L_{p_\alpha}}^\gamma \right) \|f_s\|_{L_{p_s}} \|f_t\|_{L_{p_t}} \notag \\
    &\hspace{+20mm}=\log(e+|y|)^{\sum_{k\in \{1,\ldots,n+1\}\smallsetminus \{s,t\}} \frac{1}{p_k}} \prod_{k=1}^{n+1} \|f_k\|_{L_{p_k}} . \label{eq:finalStep}
\end{align}
Behold that the exponent on the logarithm is exactly \eqref{eq:lambdast}, which is what we needed to show.

To conclude this section, we return to the question of what happens when one of the $1/p_k$ values exceeds (or equals) $\frac12$, say $1/p_u \geqslant 1/2$ for some $u\in \{1,\ldots,n+1\}$. If $u\neq s,t$, then we can simply choose $\alpha = u$ and repeat the above steps with $J_0$ being empty. If, on the other hand, $u = s$, then we actually do not need any tricks for removing the shift. Returning to Remark \ref{rem:tau=st}, we let $t = \tau$, which leaves us with the exponent on the logarithm being\footnote{The inequality just follows because $1/p_u \geqslant 1/2$, so the sum of all of the other reciprocals does not exceed $1/2$.}
\begin{align*}
    \lambda_{(u,t)}'' = \left|\frac12 - \frac{1}{p_u} \right| + \sum_{k\neq u,t} \frac{1}{p_k} \leqslant \frac{1}{p_u} ,
\end{align*}
so this case will not be the maximum exponent \eqref{eq:lambdathm} anyway. The case of $u = t$ is similar.

\begin{rem}\label{rem:sqormaxnotand}
It is tempting to ask what would happen if we could, instead, move all of the shifts onto the square functions instead of the maximal functions. If this were possible, then we would get the logarithmic exponent
\begin{align*}
    \left| \frac12 - \frac{1}{p_t} \right| + \left| \frac12 - \frac{1}{p_s} \right| ,
\end{align*}
which, when allowed to range over $(s,t)\in\cI$, clearly has \eqref{eq:lambdathm} as a maximum. That is to say, we lose sharpness in the proof when we combine the use of shifted square and maximal functions. By using purely one or the other, we get the same, sharp, result. However, it seems to not be possible to move all of the shifts onto the square functions except in the bilinear and trilinear cases, because this change of variables trick can only be used to remove the shift from one component and no more.
\end{rem}


\subsection{Some Edge Cases}\label{sec:face}

We now repeat the steps of the previous section, making small adjustments to accommodate one of the $p$-values being infinity. We will prove the result for the point 
\begin{align*}
    \left(\frac{1}{p_1},\ldots,\frac{1}{p_\nu},\ldots,\frac{1}{p_n},\frac{1}{p_{n+1}}\right) ,
\end{align*}
where $1\leqslant \nu\leqslant n+1$ and $p_\nu = \infty$. We assume that $2<p_k<\infty$ for $k\neq \nu$, i.e. each $1/p_k < 1/2$; the case, where there exists a $u$ with $1/p_u\geqslant 1/2$ can be handled in the same way it was last section. There are two situations to consider: when $\nu$ is $s$ or $t$ and when it is not. First, we will look at the latter.

As before, to be concrete, we take $\tau = t$, noting that an analogous proof works if we, instead, took $\tau = s$. We will reuse the notation of last section, adapting it slightly: define
\begin{align*}
    J_0 \subset \{1,\ldots,n+1\}\smallsetminus \{s,t,\nu\}
\end{align*}
such that $\sum_{k\in J_0\cup \{s\}} \frac{1}{p_k} < \frac12$ and there exists $\alpha \in (J_0)^c$ such that $\sum_{k\in J_0\cup \{s,\alpha\}} \frac{1}{p_k} > \frac12$.\footnote{The procedure for choosing such a set can work in the same way as we mentioned in Section \ref{sec:center}. Not including $\nu$ in $J_0$ is no issue, because $\frac{1}{p_\nu} = 0$ anyway.} It is, again, a possibility that the $p$-values are such that we can choose a $J_0$ with the sum being exactly $\frac12$; we handle this in the same way as explained in the last section. As before, we decompose
\begin{align*}
    \frac{1}{p_\alpha} = \frac{1}{q_0} + \frac{1}{q_1}
\end{align*}
so that
\begin{align*}
    \left( \sum_{k\in J_0\cup \{s\}} \frac{1}{p_k} \right) + \frac{1}{q_0} = \frac12 .
\end{align*}
Additionally, we let $0<\gamma<1$ be defined through $q_0 \gamma = p_\alpha$, implying $q_1 (1-\gamma) = p_\alpha$. We also write $J_1 = (J_0)^c\smallsetminus \{\alpha\}$.

To proceed, we have to repeat some of the steps made in Section \ref{sec:vertex}. By Lemma \ref{lem:shiftedMaxPeetre} applied to $f_\nu$, for $\sigma > d$, there exists a sequence $(S_Q)_{Q\in\cD}$ with $S_Q\subset Q$ for all $Q\in\cD$ with $|S_Q|\geqslant \frac12 |Q|$ and
\begin{align}
    \left\| \left\{ \sum_{Q\in \cD_\ell} \left( \inf_{z\in Q} \Rfk_{\sigma,\ell}^{y_\nu} f_\nu(z) \right) 1_{S_Q} \right\}_{\ell\in\nZ} \right\|_{L_\infty(\ell_\infty)} \lesssim \|f_\nu\|_{L_\infty} , \label{eq:fnu1}
\end{align}
where we have defined $\Rfk_{\sigma,\ell}^{y_\nu}$ in \eqref{eq:QRfk}.

We can bound the left hand side of \eqref{eq:logEstimate} by
\begin{align}
    &\sum_{\ell\in\nZ} \sum_{Q\in\cD_\ell} \int_Q \left( \prod_{k\neq s,t} \Rfk_{\sigma,\ell}^{y_k} f_k(x) \right) \Qfk_{\sigma,\ell}^{y_s} f_s(x) \Qfk_{\sigma,\ell} f_t(x) \, dx \notag \\
    &\hspace{+5mm}\lesssim \sum_{\ell\in\nZ} \sum_{Q\in\cD_\ell} \left( \prod_{k\neq s,t} \left( \inf_{z\in Q} \Rfk_{\sigma,\ell}^{y_k} f_k (z) \right) \right) \left( \inf_{z\in Q} \Qfk_{\sigma,\ell}^{y_s} f_s(z) \right) \left( \inf_{z\in Q} \Qfk_{\sigma,\ell} f_t(z) \right) |Q| \notag \\
    &\hspace{+5mm}\leqslant \sum_{\ell\in\nZ} \int_{\nR^d} \left( \prod_{k\neq s,t,\nu} \Rfk_{\sigma,\ell}^{y_k} f_k(x) \right) \Qfk_{\sigma,\ell}^{y_s} f_s(x) \Qfk_{\sigma,\ell} f_t(x) \left( \sum_{Q\in\cD_\ell} \left( \inf_{z\in Q} \Rfk_{\sigma,\ell}^{y_\nu} f_\nu(z) \right) 1_{S_Q}(x) \right) \, dx \label{eq:facemid} .
\end{align}
We have used \eqref{eq:PeetreisConstant} above in the second line. Note, also, that we do need to use this argument with the Peetre maximal function on each function, not just $f_\nu$, because of the need to apply \eqref{eq:PeetreisConstant}.

From this point, the process from the previous section of moving the shift off of the remaining square function is unchanged. The only difference now is that we have an extra Peetre maximal function, but it is easy to see
\begin{align*}
    \Qfk_{\sigma,\ell}^{y_s} f_s(x) = \Qfk_{\sigma,\ell} f_s(x-2^{-\ell}y_s)
\end{align*}
and similarly for the $\Rfk_{\sigma,\ell}^{y_k} f_k(x)$. In other words, Lemma \ref{lem:changeVariables} applies just as well. The only other change required in this proof is the additional use of \eqref{eq:FeffermanSteinPeetre} to get rid of the Peetre maximal functions once we have applied this change of variables trick and separated each part into its corresponding $L_{p_k}$ norm. 

Marching onward, we estimate \eqref{eq:facemid} with
\begin{align*}
    &\sum_{\ell\in\nZ} \int_{\nR^d} \left( \Qfk_{\sigma,\ell}^{y_s} f_s(x) \Rfk_{\sigma,\ell}^{y_\alpha} f_\alpha(x)^\gamma \prod_{k\in J_0} \Rfk_{\sigma,\ell}^{y_k} f_k(x) \right) \Rfk_{\sigma,\ell}^{y_\alpha} f_\alpha(x)^{1-\gamma} \left( \prod_{k\in J_1} \Rfk_{\sigma,\ell}^{y_k} f_k(x) \right) \\
    &\hspace{+25mm}\cdot\Qfk_{\sigma,\ell} f_t(x) \left( \sum_{Q\in\cD_\ell} \left( \inf_{z\in Q} \Rfk_{\sigma,\ell}^{y_\nu} f_\nu(z) \right) 1_{S_Q}(x) \right) \, dx \\
    &\hspace{+5mm}\leqslant \left\| (\Qfk_{\sigma,\ell} f_s) (\Rfk_{\sigma,\ell}^{y_\alpha - y_s} f_\alpha)^\gamma \prod_{k\in J_0} \Rfk_{\sigma,\ell}^{y_k - y_s} f_k \right\|_{L_2(\ell_2)} \| \Rfk_{\sigma,\ell}^{y_\alpha} f_\alpha \|_{L_{p_\alpha}(\ell_\infty)}^{1-\gamma} \left(\prod_{k\in J_1} \|\Rfk_{\sigma,\ell}^{y_k} f_k \|_{L_{p_k}(\ell_\infty)}\right) \\
    &\hspace{+25mm}\cdot\|\Qfk_{\sigma,\ell} f_t\|_{L_{p_t}(\ell_2)} \left\| \sum_{Q\in\cD_\ell} \left( \inf_{z\in Q} \Rfk_{\sigma,\ell}^{y_\nu} f_\nu(z) \right) 1_{S_Q} \right\|_{L_{\infty}(\ell_\infty)} ,
\end{align*}
where we have applied Lemma \ref{lem:changeVariables} in the last step. A final application of H\"older's inequality yields the following bound:
\begin{align*}
    &\|\Qfk_{\sigma,\ell} f_s\|_{L_{p_s}(\ell_2)} \|\Rfk_{\sigma,\ell}^{y_\alpha - y_s} f_\alpha \|_{L_{p_\alpha}(\ell_\infty)}^\gamma \left( \prod_{k\in J_0} \|\Rfk_{\sigma,\ell}^{y_k - y_s} f_k\|_{L_{p_k}(\ell_\infty)} \right) \| \Rfk_{\sigma,\ell}^{y_\alpha} f_\alpha \|_{L_{p_\alpha}(\ell_\infty)}^{1-\gamma} \\
    &\hspace{+20mm}\cdot\left(\prod_{k\in J_1} \|\Rfk_{\sigma,\ell}^{y_k} f_k \|_{L_{p_k}(\ell_\infty)}\right) \|\Qfk_{\sigma,\ell} f_t\|_{L_{p_t}(\ell_2)} \left\| \sum_{Q\in\cD_\ell} \left( \inf_{z\in Q} \Rfk_{\sigma,\ell}^{y_\nu} f_\nu(z) \right) 1_{S_Q} \right\|_{L_{\infty}(\ell_\infty)} .
\end{align*}
From here, we apply \eqref{eq:fnu1} to the term with $f_\nu$. Also, using \eqref{eq:FeffermanSteinPeetre} to get rid of the Peetre maximal functions in each other term, followed by Lemmas \ref{lem:shiftedsquare} and \ref{lem:shiftedmax}, gives that the final term is bounded as in \eqref{eq:finalStep}, which completes the proof of this case.

We have left ourselves one loose end, which is when $\nu = s$ or $\nu = t$. We will be brief in what follows, because it is mostly repetitious of the argument we have just given. Say, for the sake of being concrete, that $\nu = t$ (an analogous argument works if $\nu = s$). In Section \ref{sec:center}, we allowed $\tau$ to arbitrarily be one of $s$ or $t$, but here, we must choose $\tau = t = \nu$. Now, the set $J_0$ is chosen as before,
\begin{align*}
    J_0\subset \{1,\ldots,n+1\}\smallsetminus \{s,t\}
\end{align*}
so that $\sum_{k\in J_0\cup \{s\}} \frac{1}{p_k} < \frac12$ and there exists $\alpha \in (J_0)^c$ such that $\sum_{k\in J_0\cup \{s,\alpha\}} \frac{1}{p_k} > \frac12$. If there is a $J_0$ so that the sum is exactly $\frac12$, we handle this as before. Further, we decompose $1/p_\alpha = 1/q_0 + 1/q_1$ with the same criterion on $q_0$, \eqref{eq:exactly1/2}, and define $0<\gamma<1$ in the same way as above.

We use Lemma \ref{lem:shiftedSquarePeetre} (without the shift) on $f_t$: for $\sigma > d$, there exists a sequence $(S_Q)_{Q\in\cD}$ with $S_Q\subset Q$ for all $Q\in\cD$ with $|S_Q|\geqslant \frac12 |Q|$ and
\begin{align}
    \left\| \left\{ \sum_{Q\in \cD_\ell} \left( \inf_{z\in Q} \Qfk_{\sigma,\ell} f_t(z) \right) 1_{S_Q} \right\}_{\ell\in\nZ} \right\|_{L_\infty(\ell_2)} \lesssim \|f_t\|_{L_\infty} , \label{eq:fnu2}
\end{align}
where we have defined $\Qfk_{\sigma,\ell}$ in \eqref{eq:QRfk}. With a similar procedure as above, we bound \eqref{eq:logEstimate} with
\begin{align*}
    \sum_{\ell\in\nZ} \int_{\nR^d} \left( \prod_{k\neq s,t} \Rfk_{\sigma,\ell}^{y_k} f_k(x) \right) \Qfk_{\sigma,\ell}^{y_s} f_s(x) \left( \sum_{Q\in\cD_\ell} \left( \inf_{z\in Q} \Qfk_{\sigma,\ell} f_t(z) \right) 1_{S_Q}(x) \right) \, dx .
\end{align*}
Grouping terms then applying H\"older's inequality and Lemma \ref{lem:changeVariables}, this is bounded by
\begin{align*}
    &\left\| (\Qfk_{\sigma,\ell} f_s) (\Rfk_{\sigma,\ell}^{y_\alpha - y_s} f_\alpha)^\gamma \prod_{k\in J_0} \Rfk_{\sigma,\ell}^{y_k - y_s} f_k \right\|_{L_2(\ell_2)} \| \Rfk_{\sigma,\ell}^{y_\alpha} f_\alpha \|_{L_{p_\alpha}(\ell_\infty)}^{1-\gamma} \left(\prod_{k\in J_1} \|\Rfk_{\sigma,\ell}^{y_k} f_k \|_{L_{p_k}(\ell_\infty)}\right) \\
    &\hspace{+25mm}\cdot \left\| \sum_{Q\in\cD_\ell} \left( \inf_{z\in Q} \Qfk_{\sigma,\ell} f_t(z) \right) 1_{S_Q} \right\|_{L_{\infty}(\ell_2)} \\
    &\hspace{+5mm}\leqslant \|\Qfk_{\sigma,\ell} f_s\|_{L_{p_s}(\ell_2)} \|\Rfk_{\sigma,\ell}^{y_\alpha - y_s} f_\alpha \|_{L_{p_\alpha}(\ell_\infty)}^\gamma \left( \prod_{k\in J_0} \|\Rfk_{\sigma,\ell}^{y_k - y_s} f_k\|_{L_{p_k}(\ell_\infty)} \right) \| \Rfk_{\sigma,\ell}^{y_\alpha} f_\alpha \|_{L_{p_\alpha}(\ell_\infty)}^{1-\gamma} \\
    &\hspace{+25mm}\cdot\left(\prod_{k\in J_1} \|\Rfk_{\sigma,\ell}^{y_k} f_k \|_{L_{p_k}(\ell_\infty)}\right) \left\| \sum_{Q\in\cD_\ell} \left( \inf_{z\in Q} \Qfk_{\sigma,\ell} f_t(z) \right) 1_{S_Q} \right\|_{L_{\infty}(\ell_2)} .
\end{align*}
Applying \eqref{eq:fnu2} for the term with $f_t$ and \eqref{eq:FeffermanSteinPeetre} followed by Lemmas \ref{lem:shiftedsquare} and \ref{lem:shiftedmax} for the others, we arrive, for a final time, at \eqref{eq:finalStep}, which completes the proof.

\begin{rem}\label{rem:no2Linfty}
It is important to note that, while in Section \ref{sec:center} we were arbitrarily saying $\tau$ would be one of $s,t$ with either choice being as good as the other, in the latter part of this section, we explicitly choose $\tau = \nu = t$. This removes the need to take the shift off of the square function corresponding to the $L_\infty$ function, which is no longer possible after the manipulations we have done (essentially because of the infimum over $z\in Q$). For similar reasons, we excluded $\nu$ from $J_0$ at the beginning of the section (the complicated expression with $f_\nu$ is not compatible with Lemma \ref{lem:changeVariables}). This also clarifies why the argument breaks down when we have two or more components being taken from $L_\infty$.
\end{rem}


\subsection{Interpolating for the Remaining Edge Cases}\label{sec:interpolation}


We will make light use of standard terminology/notation surrounding the complex method of interpolation, see \cite[Chapter 4]{bergh2012interpolation}. The relevant compatible couple of Banach spaces will be $H_p$ and $H_q$. It was proven, see \cite[Section 3]{calderon1977parabolic} (for $q<\infty$ below) and \cite[Theorem 1]{janson1982interpolation} (for $q=\infty$), that
\begin{align*}
    (H_p,H_q)_\theta = H_r , \quad \frac{1}{r} = \frac{1-\theta}{p} + \frac{\theta}{q} ,
\end{align*}
for $0\leqslant \theta\leqslant 1$ and $0<p<q\leqslant\infty$. Corresponding to this, we have the following multilinear interpolation theorem, which is just a special case of \cite[Theorem 4.4.1]{bergh2012interpolation}.

\begin{prop}\label{prop:multiInterp}
Consider two sets of $p$-values: $1\leqslant p_{1,j},\ldots,p_{n+1,j} \leqslant \infty$ for $j=0,1$. Suppose
\begin{align*}
    \Gamma : \bigoplus_{k=1}^{n+1} (H_{p_{k,0}}\cap H_{p_{k,1}})\to \nC
\end{align*}
is multilinear and there exist $M_0,M_1 > 0$ such that $\Gamma$ satisfies
\begin{align*}
    |\Gamma(f_1,\ldots,f_{n+1})| \leqslant M_j \prod_{k=1}^{n+1} \|f_k\|_{H_{p_{k,j}}} , \quad f_k\in H_{p_{k,0}}\cap H_{p_{k,1}} ,
\end{align*}
for $j=0,1$. Then, for every $0\leqslant\theta\leqslant 1$, there is a unique extension of $\Gamma$ to a multilinear operator,
\begin{align*}
    \Gamma : \bigoplus_{k=1}^{n+1} (H_{p_{k,0}},H_{p_{k,1}})_\theta \to \nC ,
\end{align*}
with operator norm at most $M_0^{1-\theta} M_1^\theta$.
\end{prop}

Before we interpolate we need to make two remarks on our process. First, as was the case in the last few sections, we need only establish \eqref{eq:logEstimate}. For this reason, we define $\Gamma_{(s,t)}$ to be the left hand side of \eqref{eq:logEstimate} (without the absolute values) for $(s,t)\in\cI$. By interpolating with $\Gamma_{(s,t)}$ we have the logarithms as our operator norms, which will be essential. For the rest of this subsection, $(s,t)\in\cI$ will be fixed and we will just write $\Gamma$ instead of $\Gamma_{(s,t)}$.

Second, the reader may have noticed that we have so far only worked with $\Gamma$ for functions in $\cS_0(\nR^d)$. Therefore, in order to apply the interpolation result, we must first extend the vertex cases from Section \ref{sec:vertex}. Explicitly, let $j\in \{1,\ldots,n+1\}$ and say that $p_k = \infty$ for $k\neq j$ and $p_j = 1$; using the fact that
\begin{align*}
    |\Gamma(f_1,\ldots,f_{n+1})| \lesssim \log(e + |y|) \|f_j\|_{H_1} \prod_{k\neq j} \|f_k\|_{L_\infty}
\end{align*}
for $f_1,\ldots,f_{n+1}\in\cS_0(\nR^d)$, we extend $\Gamma$ so that the above holds for $f_k\in L_\infty$ for $k\neq j$ and $f_j\in H_1$. If we did not have to deal with $L_\infty$, this would be a trivial extension by density, but $\cS_0$ is not dense in $L_\infty$, so a more complicated, but still standard, process (using the transpose operators) must be undertaken. Since it is standard, we omit this and invite the reader to check \cite[Section 4.4.1]{muscalu2013classical}.

We can now proceed with the interpolation: the idea is to iterate Proposition \ref{prop:multiInterp}---first, with the points $(1,0,\ldots,0)$ and $(0,1,0,\ldots,0)$, which establishes the sharp result for points of the form $(\frac{1}{p_1},\frac{1}{p_2},0,\ldots,0)$, second with the arbitrary point $(\frac{1}{p_1},\frac{1}{p_2},0,\ldots,0)$ and $(0,0,1,0,\ldots,0)$, which gives the sharp result for points with all but the first three components being zero, and so on. For simplicity, we will establish the result for points of the form $(\frac{1}{p_1},\ldots,\frac{1}{p_j},0,\ldots,0)$ with $2\leqslant j\leqslant n-1$;\footnote{We note that these calculations are void when $n=2$, but this section is also unnecessary when $n=2$, so we can just take $n>2$.} those cases with zeros in other locations follow similarly.

We are lucky here in two ways: one, the points we are using for interpolation are simple to work with, and two, the operator norm in each case will just be bounded by $\log(e+|y|)$, meaning that the logarithmically convex combination---$M_0^{1-\theta} M_1^{\theta}$---from Proposition \ref{prop:multiInterp} will just be $\log(e+|y|)$ each time we interpolate.

From Section \ref{sec:vertex} (and the extension discussed above), we know that
\begin{align*}
    \Gamma : H_1 \times L_\infty \times \cdots \times L_\infty \to \nC , \quad \Gamma : L_\infty \times H_1 \times L_\infty\times \cdots \times L_\infty \to \nC ,
\end{align*}
both with operator norms bounded by a constant multiple of $\log(e+|y|)$. Therefore, we have taken, using the notation of Proposition \ref{prop:multiInterp}, $p_{1,0} = p_{2,1} = 1$ and the rest to be infinity.

Performing the first interpolation, we see that $\Gamma : L_{p_1} \times L_{p_2} \times L_\infty \times \cdots \times L_\infty \to \nC$ with operator norm bounded by a constant multiple of $\log(e+|y|)$.\footnote{We have used, additionally, that if $A_0,A_1$ are a compatible Banach couple, then $(A_0,A_1)_\theta = (A_1,A_0)_{1-\theta}$ for $0\leqslant \theta\leqslant 1$. Check \cite[Theorem 4.2.1]{bergh2012interpolation} for a proof of this simple result.} Thus, we have proven the sharp result in the case where all but the first two of the components are taken from $L_\infty$.

If $n=3$, we are already done, so now we proceed with $n > 3$. Let us now assume that, for some $j\in \{2,\ldots,n-2\}$ and any $1<p_1,\ldots,p_j<\infty$,\footnote{If we wanted to take one of the $p$-values to be $1$, then that would be a vertex point, which is known. If we wanted more than $(n+1)-j$ of them to be $\infty$, that would be covered by earlier steps in the induction.}
\begin{align*}
    \begin{cases}
        \Gamma : L_{p_1}\times \cdots \times L_{p_j} \times L_\infty\times \cdots \times L_\infty\to \nC , \\
        \Gamma : L_\infty \times \cdots \times L_\infty\times H_1 \times L_\infty\times \cdots \times L_\infty \to \nC ,
    \end{cases}
\end{align*}
where $H_1$ is in the $(j+1)$-st spot. In the latter case, we know from Section \ref{sec:vertex} that the operator norm is at most a constant multiple of $\log(e+|y|)$; in the former, we are assuming the same. Interpolating yields, then, that
\begin{align*}
    \Gamma : L_{q_1}\times \cdots \times L_{q_{j+1}} \times L_\infty \times \cdots \times L_\infty \to \nC ,
\end{align*}
with operator norm bounded, up to a constant, by $\log(e+|y|)$ and where
\begin{align*}
    \frac{1}{q_k} = \frac{1-\theta}{p_k} , \quad \frac{1}{q_{j+1}} = 1-\theta
\end{align*}
for $k=1,\ldots,j$ and a $0<\theta<1$ of our choice. Clearly, since $1<p_1,\ldots,p_j<\infty$ and $0<\theta<1$ can be chosen at will, we have the sharp result for the case where all but $j+1$ of the components are taken from $L_\infty$. This finishes the proof.

\begin{rem}\label{rem:interpLimit}
Observe that if we try the above process one more time to get the result for the points from Section \ref{sec:face}, i.e. those of the form $( \frac{1}{p_1} , \ldots , \frac{1}{p_n} , 0 )$, for $1<p_1,\ldots,p_n<\infty$, then we would get an operator norm bounded by $\log(e + |y|)$, which no longer has the sharp exponent, \eqref{eq:lambdathm}. More generally, the interpolation yields the sharp result exactly when we interpolate two points for which the largest $n-1$ of their coordinates align.\footnote{For example, if $n=3$, we could interpolate $(\frac12,\frac14,\frac14,0)$ with $(\frac14,\frac12,0,\frac14)$ and get the sharp result, but not $(\frac12,\frac14,\frac14,0)$ with $(\frac12,0,\frac16,\frac13)$.} This is due to the logarithmic convexity of the operator norm when applying complex interpolation combined with the structure of \eqref{eq:lambdathm}.
\end{rem}

\begin{rem}
It is also possible to get the full result via interpolation by continuing from what we have shown above, also interpolating with the two points
\begin{align*}
    \left( \frac{1}{n+1} , \ldots , \frac{1}{n+1} \right) , \quad \left( \frac{1}{n} ,\ldots,\frac{1}{n} , 0 \right)
\end{align*}
(or with the zero in the latter point in different places). In fact, one can easily see based on Remark \ref{rem:interpLimit}, that these are the only two points one can use if one wishes to prove the full sharp result. Further, the sharp result cannot be proven for these points based on interpolation with other points because of what we discussed in Remark \ref{rem:interpLimit}.
\end{rem}

%% file: Sections/counterexample.tex
\section{Proof of Sharpness}\label{sec:converse}

In this final section, we give a counterexample, which establishes the converse direction stated in Theorem \ref{thm:main}.


\subsection{The Counterexample}\label{sec:counterexample}

It turns out that a simple adaptation of the counterexample from \cite{dosidis2024boundedness}, which proved the sharpness of the bilinear result, works again in the multilinear setting. We will be brief, because much of what is below is repetitious of \cite[Section 5.2]{dosidis2024boundedness}.

Throughout this section, we will be looking at the $n$-linear operator, $T$, rather than its $(n+1)$-linear counterpart, $\Lambda$. For that reason, when we write, e.g., $j\neq s,t$, we mean $j\in \{1,\ldots,n\}\smallsetminus \{s,t\}$, i.e. $j=n+1$ is not a possibility.

Let $\zeta_k = 10k$ for $k\in\nN$ and write $e_1$ to mean the standard unit vector $(1,0,\ldots,0)\in \nR^d$. Suppose, to be suggestive, that the smallest two of the reciprocal $p$-values are $1/p_s$ and $1/p_t$ for some $1\leqslant s < t \leqslant n$. That is, for now, $1/p'$ is assumed to not be one of the two smallest reciprocal $p$-values; we will return to this point at the end of the section.

We model $K$ off of one of the pieces into which we decomposed it in \eqref{eq:Kdecomp}, so for that we define $\beta,\eta\in\cS$ to be radial functions satisfying the following properties:
\begin{align*}
    \supp(\widehat{\eta}) \subset \left\{\xi\in\nR^d : |\xi| \leqslant \frac{1}{100} \right\} , \quad \widehat{\eta}(\xi) = 1 \, \text{ on } \, \left\{\xi\in\nR^d : |\xi|\leqslant \frac{1}{200} \right\}
\end{align*}
and
\begin{align*}
    \supp(\widehat{\beta}) \subset \left\{\xi\in\nR^d : \frac{10}{11} \leqslant |\xi|\leqslant \frac{11}{10} \right\} , \quad \widehat{\beta}(\xi) = 1 \, \text{ on } \, \left\{\xi\in\nR^d : \frac{20}{21}\leqslant |\xi|\leqslant \frac{21}{20} \right\} .
\end{align*}
Due to these properties, we observe that
\begin{align}
    \widehat{\beta}(2^{-\ell}\xi) \widehat{\eta}(\xi - 2^{\zeta_k}e_1) = \begin{cases}
        0 & \ell\neq \zeta_k , \\
        \widehat{\eta}(\xi - 2^{\zeta_k}e_1) & \ell = \zeta_k .
    \end{cases} \label{eq:betaEtaProduct}
\end{align}
We now build $K,f_1,\ldots,f_n$ with these functions. Define
\begin{align*}
    f_s(x) = \sum_{k=1}^N \eta(x + 2^{\zeta_N - \zeta_k}e_1) e^{2\pi i(x,2^{\zeta_k}e_1)} , \quad f_t(x) = \sum_{k=1}^N \eta(x + 2^{\zeta_N - \zeta_k}e_1) e^{-2\pi i(x,2^{\zeta_k}e_1)} ,
\end{align*}
and $f_j(x) = \beta(x)$ for $j\neq s,t$. One can check that each of these has Fourier transform supported away from the origin, so they are in $\cS_0(\nR^d)$. The kernel is
\begin{align*}
    K(y_1,\ldots,y_n) = \beta(y_s - 2^{\zeta_N}e_1)\beta(y_t - 2^{\zeta_N}e_1) \prod_{j\neq s,t} \eta(y_j) .
\end{align*}
Obviously, $\|f_j\|_{H_{p_j}}\lesssim 1$ for $j\neq s,t$. Following \cite{dosidis2024boundedness}, we have $D_\lambda(K)\lesssim N^\lambda$ and
\begin{align*}
    \|f_s\|_{H_{p_s}} \lesssim N^{1/p_s} , \quad \|f_t\|_{H_{p_t}} \lesssim N^{1/p_t} ,
\end{align*}
using the characterization of $L_p$ through square functions, see \cite[Equation (4.14)]{dosidis2024boundedness}.

It remains to show $\|T(f_1,\ldots,f_n)\|_{Y_p} \gtrsim N$. To do so, we take the Fourier transform of
\begin{align*}
    \sum_{\ell\in\nZ} K_\ell * (f_1\otimes\cdots\otimes f_n)(z_1,\ldots,z_n) ,
\end{align*}
witness an amazing cancellation based on \eqref{eq:betaEtaProduct}, take the inverse Fourier transform, evaluate it at $(x,\ldots,x)$, and discover that
\begin{align*}
    T(f_1,\ldots,f_n)(x) = N(\eta(x))^2 \beta(x)^{n-2} .
\end{align*}
The calculation is almost identical to the one in \cite{dosidis2024boundedness}, so we skip the details.

From here, it is obvious that $\|T(f_1,\ldots,f_n)\|_{Y_p} \gtrsim N$. Pulling this all together,
\begin{align*}
    \frac{\|T(f_1,\ldots,f_n)\|_{Y_p}}{D_\lambda(K) \|f_1\|_{H_{p_1}} \cdots \|f_n\|_{H_{p_n}}} \gtrsim N^{1 - \left(\frac{1}{p_s} +\frac{1}{p_t} \right) - \lambda} = N^{\left(\sum_{j\neq s,t} \frac{1}{p_j} \right) - \lambda} ,
\end{align*}
which establishes our claim.

\begin{rem}
It may seem somehow unsatisfying that most of the functions are not doing any work for us (by which we mean that we only have $\|f_j\|_{H_{p_j}}\lesssim 1$ for $j\neq s,t$). This is, however, a reflection of the proof of the positive direction of the theorem. Imagining $N$ as a stand-in for $\log(e+|y|)$, we see that the construction of the counterexample above has brought this factor exactly out of the two functions attached to the two smallest reciprocal $p$-values, leaving the others untouched. This corresponds to using only shifted square function estimates, the possibility of which we mentioned in Remark \ref{rem:sqormaxnotand}.

We could ask, from this perspective, if there is a different counterexample, which mirrors the fact that in our proof we have shifts only on the maximal functions. To that question, we see it as likely that such a counterexample exists, but we do not pursue a construction here.
\end{rem}

The final issue to address is what happens if $1/p'$ is one of the two smallest reciprocal $p$-values; say that the other is $1/p_s$ for some $1\leqslant s\leqslant n$. Returning to the adjoint operators, proving the positive result,
\begin{align}
    \|T(f_1,\ldots,f_n)\|_{Y_p} \leqslant C D_\lambda(K) \|f_1\|_{H_{p_1}} \cdots \|f_n\|_{H_{p_n}} , \label{eq:normalpositivedirec}
\end{align}
is, of course, equivalent to proving
\begin{align}
    &\|T^j(f_1,\ldots,f_{j-1},f_{n+1},f_{j+1},\ldots,f_n)\|_{Y_p} \notag \\
    &\hspace{+10mm}\leqslant C D_\lambda(K) \|f_1\|_{H_{p_k}}\cdots \|f_{j-1}\|_{H_{p_{j-1}}} \|f_{n+1}\|_{H_{p'}} \|f_{j+1}\|_{H_{p_{j+1}}} \cdots \|f_n\|_{H_{p_n}} , \label{eq:counterwhent=n+1}
\end{align}
for any $1\leqslant j\leqslant n$ (with obvious modifications to the above if $j=1$ or $j=n$). Therefore, we can just choose any $j\neq s$ and then reproduce our counterexample, but for \eqref{eq:counterwhent=n+1} rather than \eqref{eq:normalpositivedirec}. With that, we are done.